\documentclass[twoside,12pt]{article}
\usepackage{amsmath}
\usepackage{amssymb}
\usepackage{cite}
\usepackage[usenames,dvipsnames]{color}
\definecolor{rosso}{rgb}{0.8,0,0}

\def\pier{\color{rosso}}
\def\revis{\color{red}}
\def\colli{\color{blue}}

\let\pier\relax
\let\revis\relax
\let\colli\relax

\title{{\pier Equation and dynamic boundary condition}\\
of {C}ahn--{H}illiard type with singular potentials}

\author{Pierluigi Colli\\
Dipartimento di Matematica, Universit\`a degli Studi di Pavia\\
Via Ferrata~1, 27100 Pavia, Italy\\
E-mail: \texttt{pierluigi.colli@unipv.it}\\
\and \\ Takeshi Fukao\\
Department of Mathematics, Faculty of Education\\
Kyoto University of Education\\
1~Fujinomori, Fukakusa, Fushimi-ku, Kyoto~612-8522 Japan\\
E-mail: \texttt{fukao@kyokyo-u.ac.jp}}
 
\date{}

\newcommand\testopari{\sc {\pier Pierluigi Colli and Takeshi Fukao}}
\newcommand\testodispari{\sc {\pier Equation and boundary condition
of {C}ahn--{H}illiard type}}
\begin{document}

\maketitle

\begin{abstract}
The well-posedness of a system of partial differential equations and dynamic 
boundary conditions, both of {C}ahn--{H}illiard type, is discussed. 
The existence of a weak solution and its continuous dependence on the data are 
proved using {\pier a} suitable setting for the conservation of a total mass in the bulk {\pier plus} 
the boundary. A very general class of double-well like potentials is allowed. 
Moreover, some further regularity is obtained to guarantee the strong solution.

\vspace{2mm}
\noindent \textbf{Key words:}~~{C}ahn--{H}illiard {\pier system}, dynamic boundary condition, mass conservation, {\pier well-posedness, strong solution.} 

\vspace{2mm}
\noindent \textbf{AMS (MOS) subject clas\-si\-fi\-ca\-tion:} {\pier 35K61, 35K25, 35D30, 35D35,} 80A22.

\end{abstract}

\section{Introduction}
\setcounter{equation}{0}

The {C}ahn--{H}illiard equation \cite{CH58, EZ86} {\pier yields} a 
famous description of the evolution phenomena on the solid-solid phase separation. 
In general, {\pier an} evolution process goes on diffusively. However, the {\pier phenomenon of the} 
solid-solid phase separation does not seem to follow on this structure{\pier :}
more precisely, each phase {\pier concentrates and this process is usually known as spinodal} decomposition. 
The {C}ahn--{H}illiard equation is a {\pier celebrated} model which describes 
this decomposition by the simple framework of partial differential equations. 
Thereon, the volume conservation is {\pier a key property of the structure: you
can observe the pattern formation that is restricted by the property of conservation on the decomposition.} On the other hand, in the real world 
there are many phenomena of pattern formation which do not have the structure of conservation. 
However, there is a possibility that actually the structure of conservation 
{\pier is hidden somewhere, and we only cannot find it at} the level of observation. 

In this paper, {\pier a system coupling the same kind of equations and dynamic 
boundary conditions of {C}ahn--{H}illiard type is investigated. We aim to describe it at 
once. Let $0<T<+\infty$ be some fixed time and let $\Omega \subset \mathbb{R}^{d}$, 
$d=2$ or $3$, be a bounded smooth domain occupied by a material: the boundary $
\Gamma$ of $\Omega $ is supposed to be smooth enough as well.  
We start from the following equations of 
{C}ahn--{H}illiard type in the domain $Q:=\Omega \times (0,T)$
\begin{gather} 
	{\pier \partial_t u} -\Delta \mu = 0 
	\quad \mbox{in }Q, 
	\label{CH1}
	\\
	\mu = -\Delta u + W'(u)-f 
	\quad \mbox{in }Q,
	\label{CH2}
\end{gather}
{\pier where ${\pier \partial_t}$ denotes the partial derivative with respect to time and, as usual, $\Delta$ represents  the Laplacian operator acting on the space variables. Here, the
unknowns $u $ and $ \mu:Q \to \mathbb{R}$ 
stand for the order parameter and} 
the chemical potential, respectively. 
In order to consider the dynamics on the boundary $\Sigma:=\Gamma\times (0,T)$, 
we also introduce the unknowns 
$u_\Gamma, \mu_\Gamma:\Sigma \to \mathbb{R}$
such that
\begin{gather} 
	u_\Gamma =u _{|_\Gamma }, \quad \mu _\Gamma =\mu _{|_\Gamma }
	\quad \mbox{on }\Sigma, 
	\label{trace}
\end{gather} 
$u_{|_\Gamma }$ and $\mu_{|_\Gamma }$ being the traces of $u$ and $\mu$,} 
and consider the same type of equations on the boundary
\begin{gather} 
	{\pier \partial_t}  u_\Gamma +\partial_\nu \mu -\Delta _\Gamma \mu _\Gamma =0
	\quad \mbox{on }\Sigma,
	\label{CHB1}
	\\
	\mu _\Gamma =\partial _\nu u - \Delta _\Gamma u_\Gamma +W'_\Gamma (u_\Gamma )-f_\Gamma 
	\quad \mbox{on }\Sigma,
	\label{CHB2}
\end{gather}
where {\pier the extra terms in \eqref{CHB1} and \eqref{CHB2} contain 
the outward normal derivative $\partial_\nu $  on $\Gamma $, and where 
$\Delta _{\Gamma }$ denotes} the {L}aplace--{B}eltrami operator
on $\Gamma $ (see, e.g., \cite[Chapter~3]{Gri09}).
We can say that this kind of dynamic boundary condition \eqref{CHB1}--\eqref{CHB2} is 
a sort of transmission problem between the 
dynamics in the bulk $\Omega $ and {\pier the} one on the boundary $\Gamma$. 
Together with the conditions
\begin{gather} 
	u(0)=u_0
	\quad 
	\mbox{in }\Omega, 
	\label{IC1}
	\\
	u_\Gamma  (0)=u_{0\Gamma}
	\quad \mbox{on }\Gamma,
	\label{IC2}
\end{gather} 
the {\pier initial and boundary value problem expressed in 
\eqref{CH1}--\eqref{IC2} is termed (P)}. 
In \eqref{CH2} and \eqref{CHB2}, the nonlinear {\pier terms} $W'$ and 
$W_\Gamma '$ play {\pier some important role, since they are the derivatives 
of the functions $W$ and $W_\Gamma $ usually referred as double-well potentials,}
with two minima and a local unstable maximum in between. 
The prototype model 
is provided by $W(r)=W_\Gamma (r)=(1/4)(r^2-1)^2$ 
so that $W'(r)=W_\Gamma '(r)=r^3-r${\pier , $r\in \mathbb{R}$,} is the sum of 
an increasing function with a power growth and another smooth 
(in particular, Lipschitz continuous)
function which breaks the monotonicity properties of 
the former and is related to the non-convex part 
of the potential $W$ or $W_\Gamma$. 
{\pier To our knowledge, this problem {\rm {\pier (P)}} was formulated  
by {G}oldstein, {M}iranville and {S}chimperna~\cite{GMS11} 
and analyzed from various viewpoints (see~\cite{CGM13, CMZ11, CP14, GM13}).  
In this paper, we treat more general cases
for such nonlinearities, that is, we let $W'$ and $W_\Gamma '$ be the sum
of a maximal monotone graph and of a Lipschitz perturbation and we are able 
to show the existence of strong solutions for our system under appropriate assumptions. Our treatment is related to the approach followed in \cite{CC13, CGS14} for some other class of problems.}

{\pier As is well known, for the usual {C}ahn--{H}illiard system~\eqref{CH1}--\eqref{CH2}, the conservation of (the mean value of) $u$ is guaranteed under the homogeneous {N}eumann boundary condition 
\begin{gather*}
	\partial_\nu \mu  = 0 
	\quad \mbox{on }\Sigma,
\end{gather*}
for $\mu $; namely, thanks to this, by simply integrating \eqref{CH1} over $\Omega$  we easily obtain from \eqref{IC1} that 
\begin{gather*}
	 \int_{\Omega }^{} u(t) dx = \int_{\Omega }^{} u_0 dx 
	\quad \mbox{for all } t \in [0,T].
\end{gather*}
Let us mention the related papers~{\revis\cite{Aik95, CGM13, KN96, Kub12}} 
on which the structure of 
conservation was treated in an abstract framework. The new issue of the 
problem {\rm {\pier (P)}} is the natural consequence of 
a mass constraint involving the values of $u$ both in the bulk and 
on the boundary. In fact, {\pier it arises as an outcome of 
\eqref{CH1} and \eqref{CHB1}} that the solution $u$ {\pier satisfies}
\begin{gather}
	\int_{\Omega }^{}u(t) dx+ \int_{\Gamma }^{} u_{\Gamma }(t) d\Gamma 
	= \int_{\Omega }^{}u_0 dx+ \int_{\Gamma }^{} u_{0\Gamma } d\Gamma
	\quad \mbox{for all }t \in [0,T].
    \label{pier0} 
\end{gather}
Concerning the model, let us point out that, 
under suitable choices of $W'$ and $W_\Gamma '$, we could flexibly 
realize some dynamics for pattern formation respecting \eqref{pier0}. 
For example, we may examine the case when every pattern is going 
to disappear, namely the bulk comes to be occupied by a single phase 
except near the boundary (occupied by another phase) and thus 
\eqref{pier0} plays as a conservation law. We invite the reader to compare
the approach of this paper with the one adopted in \cite{CF14}, where
the {\revis Allen--Cahn} equation, coupled with dynamic boundary conditions,
is investigated under a mass constraint which involves the solution inside the domain 
and its trace on the boundary. In that case, the constraint is rather imposed 
(in opposition with \eqref{pier0}, which is a gift of the problem)
and the system of nonlinear partial differential equations 
can be formulated as a variational inequality.}

A brief outline of the present paper along with 
a short description of the various items is as follows.

In Section 2, we present the main results, 
consisting in the well-posedness of the system {\pier \eqref{CH1}--\eqref{IC2}} of partial differential equations and dynamic 
boundary conditions, both of {C}ahn--{H}illiard type. 
We define a weak and strong solution of the problem {\pier (P)}. 
We also write the system as an evolution inclusion.

In Section 3, we prove the continuous dependence {\pier on the data} and
this result entails the uniqueness property. 

In Section 4, we prove the existence result. 
The proof is split in several steps. 
First, we construct an approximate solution 
by substituting the maximal monotone graphs with their 
{\pier {Y}osida regularizations.}
The solvability of the approximate problem is guaranteed by the 
abstract theory of doubly nonlinear evolution inclusions~\cite{CV90}. 
Moreover, arguing in a similar way as in \cite{CF15},
we show that the solution satisfies suitable regularity properties
{\pier and uniform estimates}. Finally, from these estimates, 
we can pass to the limit and conclude the existence proof of the weak solution. 
Next, we can proceed by considering {\pier some} additional uniform 
{\pier estimates} in order to obtain the strong solution. 

{\pier Finally, Section 5 contains an Appendix collecting some useful verifications.}

{\pier Anyway, for the reader's convenience, a detailed index of sections and subsections follows.}

\begin{itemize}
 \item[1.] Introduction
 \item[2.] Main results
\begin{itemize}
 \item[2.1.] Weak formulation
 \item[2.2.] Definition of the solution and main theorem
 \item[2.3.] Abstract formulation 
\end{itemize}
 \item[3.] Continuous dependence 
 \item[4.] Existence
\begin{itemize}
 \item[4.1.] Approximation of the problem
 \item[4.2.] A priori estimates
 \item[4.3.] Passage to the limit as $\varepsilon \to 0$
 \item[4.4.] Regularity result
\end{itemize}
 \item[5.] Appendix
\end{itemize}

\section{Main results}
\setcounter{equation}{0}

In this section, our main result is stated. 
We present our system of equations and conditions now:
\begin{gather} 
	{\partial_t u}-\Delta \mu = 0 
	\quad \mbox{a.e.\ in }Q,
	\label{GMS1}
\\
	\mu = -\Delta u + \xi + \pi (u) -f, 
	\quad \xi \in \beta (u)
	\quad \mbox{a.e.\ in }Q,
	\label{GMS2}
\\
	u_\Gamma =u _{|_\Gamma }, \quad \mu _\Gamma =\mu _{|_\Gamma },
	\quad 
	{\partial_t u_\Gamma }+\partial_\nu \mu -\Delta _\Gamma \mu _\Gamma =0
	\quad \mbox{a.e.\ on }\Sigma,
	\label{GMS3}
\\ 
	\mu _\Gamma =\partial_\nu u - \Delta _\Gamma u_\Gamma +\xi _\Gamma + \pi _\Gamma (u_\Gamma )-f_\Gamma, 
	\quad \xi_\Gamma  \in \beta _\Gamma (u_\Gamma )
	\quad \mbox{a.e.\ on }\Sigma,
	\label{GMS4}
\\
	u(0)=u_0 
	\quad 
	\mbox{a.e.\ in }\Omega, \quad 
	u_\Gamma  (0)=u_{0\Gamma}
	\quad \mbox{a.e.\ on }\Gamma,
	\label{GMS5}
\end{gather} 
where $f:Q \to \mathbb{R}$, $f_\Gamma :\Sigma \to \mathbb{R}$, 
$u_0:\Omega \to \mathbb{R}$, $u_{0\Gamma} :\Gamma \to \mathbb{R}$ are given functions; 
$\beta$ stands for the subdifferential of the convex part $\widehat{\beta }$ and 
$\pi $ stands for the derivative of the concave perturbation $\widehat{\pi}$ of a 
double well potential $W(r)=\widehat{\beta }(r)+\widehat{\pi}(r)$ for all 
$r \in \mathbb{R}$. 
Here $\beta $ is generalized to the case of 
maximal monotone graph in 
$\mathbb{R} \times \mathbb{R}$. 
$\beta _\Gamma$ and $\pi _\Gamma $ 
have the same property as $\beta $ and $\pi$, respectively. 
Typical examples of $\beta$, $\beta _\Gamma $ and $\pi $, $\pi _\Gamma $ 
are given as follows:
\begin{itemize}
	\item $\beta(r)=\beta _\Gamma (r)=r^3$, $\pi (r)=\pi _\Gamma (r)=-r$ 
	for all $r \in \mathbb{R}$ with $D(\beta )=D(\beta _\Gamma )=\mathbb{R}$ 
	for the prototype double well potential $W(r)=(r^2-1)^2/4$;
	\item $\beta(r)=\beta _\Gamma (r)=\ln((1+r)/(1-r))$, $\pi (r)=\pi _\Gamma (r)=-2cr$ 
	for all $r \in D(\beta )$ with $D(\beta )=D(\beta _\Gamma )=(-1,1)$ 
	for the logarithmic double well potential $W(r)=((1+r)\ln(1+r)+(1-r)\ln(1-r))-cr^2$
	where $c>0$ is a large constant which breaks convexity;
	\item $\beta(r)=\beta _\Gamma (r)=\partial I_{[-1,1]}(r)$, $\pi (r)=\pi _\Gamma (r)=-r$ 
	for all $r \in D(\beta )$ with $D(\beta )=D(\beta _\Gamma )=[-1,1]$ 
	for the singular potential $W(r)=I_{[-1,1]}(r)-r^2/2$.
\end{itemize}
{\pier Of course, it is not necessary that $\beta$ and $\beta_\Gamma $ are the same graph
or the same kind of graphs, what is important is that they respect the compatibility 
condition (A6),
which is stated below. Our working assumption is that the boundary potential somehow dominates 
the potential in the bulk, cf. \cite{CC13, CF14,CF15,CGS14} for analogous approaches. We also 
point out that a mixing of the first two cases is considered by the results reported in 
\cite{CGM13}: there, the inclusions in \eqref{GMS2} and \eqref{GMS4} actually reduce to the 
equalities $\xi =\beta (u)$ and $\xi _\Gamma =\beta _\Gamma (u_\Gamma )$. In the 
present paper, we can handle also the case of effective graphs.}

\subsection{Weak formulation}
We treat the problem {\rm {\pier (P)}} by a system of variational formulations. 
To this aim, we introduce {\pier the} 
spaces $H:=L^2(\Omega )$, $V:=H^1(\Omega )$, 
$H_\Gamma :=L^2(\Gamma )$, $V_\Gamma :=H^1(\Gamma )$
with usual norms $| \cdot |_{H}$, $|\cdot |_{V}$, $|\cdot |_{H_\Gamma}$, $|\cdot |_{V_\Gamma}$ 
and inner products $(\cdot,\cdot )_{H}$, $(\cdot ,\cdot )_{V}$,
$(\cdot,\cdot )_{H_\Gamma}$, $(\cdot ,\cdot )_{V_\Gamma}$, respectively. 
Moreover, {\pier we} put 
$\mbox{\boldmath $ H$}:=H \times H_\Gamma $ and 
\begin{gather*}
	\mbox{\boldmath $ V$}:=\left\{ (z,z_\Gamma ) \in V \times V_\Gamma \ : \  z_\Gamma =z_{|_\Gamma } \right\}.
\end{gather*}
Then, $\mbox{\boldmath $ H $}$ and $\mbox{\boldmath $ V $}$ are Hilbert spaces 
with the inner products 
\begin{align*}
	(\mbox{\boldmath $ u $},\mbox{\boldmath $ z $}
	)_{\mbox{\scriptsize \boldmath $ H$}}
	:=(u,z)_{H} + (u_\Gamma ,z_\Gamma )_{H_\Gamma } \quad 
	& \mbox{for all}~\mbox{\boldmath $ u$}
	:=(u,u_{\Gamma }), \,
	\mbox{\boldmath $ z$}:=(z,z_{\Gamma }) 
	\in \mbox{\boldmath $ H$},\\
	(\mbox{\boldmath $ u $},\mbox{\boldmath $ z $}
	)_{\mbox{\scriptsize \boldmath $ V$}}
	:=(u,z)_{V} + (u_\Gamma ,z_\Gamma )_{V_\Gamma } \quad 
	& \mbox{for all}~\mbox{\boldmath $ u$}
	:=(u,u_{\Gamma }), \,
	\mbox{\boldmath $ z$}:=(z,z_{\Gamma }) 
	\in \mbox{\boldmath $ V$},
\end{align*}
and related {\revis norms.} As a remark, let us restate that if 
$\mbox{\boldmath $ z$}:=(z,z_{\Gamma}) \in \mbox{\boldmath $ V$}$ then $z_{\Gamma }$ 
is exactly the trace of $z$ on $\Gamma$;
while, if $\mbox{\boldmath $ z$}:=(z,z_{\Gamma})$ is just in $ \mbox{\boldmath $ H$}$, 
then $z \in H$ and $z_{\Gamma } \in H_{\Gamma }$ are independent. 
From now on, we use the notation of a bold letter like $\mbox{\boldmath $ u$}$ to denote 
the pair which corresponds to the letter, that is $(u,u_\Gamma )$ for $\mbox{\boldmath $ u$}$. 
Now{\pier , for $\mbox{\boldmath $ z$} :=(z,z_\Gamma ) \in \mbox{\boldmath $ V$}$ and $t\in (0,T)$, we test \eqref{GMS1} by $z$ and use \eqref{GMS3} to infer}
\begin{align}
	\int_{\Omega }^{} {\partial_t u} (t)z dx
	+\int_{\Gamma }^{} {\partial_t u_\Gamma } (t)z_\Gamma  d\Gamma
	+\int_{\Omega }^{} \nabla \mu(t) \cdot \nabla z dx 
	+\int_{\Gamma }^{} \nabla _\Gamma \mu _\Gamma(t) \cdot \nabla _\Gamma z_\Gamma d\Gamma 
	=0.
	\label{wf1} 
\end{align} 
We also test \eqref{GMS2} by $z$ and {\pier exploit \eqref{GMS4};} then we obtain 
\begin{align}
	\lefteqn{  
	\int_{\Omega }^{} \mu(t) z dx + \int_{\Gamma }^{} \mu _\Gamma(t) z_\Gamma d\Gamma 
	} \nonumber \\
	& = \int_{\Omega }^{} \nabla u(t)\cdot \nabla z dx + 
	\int_{\Gamma }^{} \nabla _\Gamma u _\Gamma(t) \cdot \nabla _\Gamma z_\Gamma d\Gamma 
	 + \int_{\Omega }^{} \bigl( 
	\xi(t) +\pi \bigl( u(t) \bigr)-f(t) \bigr)z dx 
	\nonumber \\
	& \quad {} 
	+ \int_{\Gamma }^{} \bigl( 
	\xi_\Gamma(t)  +\pi_\Gamma  \bigl( u_\Gamma (t) \bigr)-f_\Gamma(t) \bigr )z_\Gamma  d\Gamma.
	\label{wf2}
\end{align}
Now, let us take $\mbox{\boldmath $ z$}=\mbox{\boldmath $ 1$}:=(1,1)$ in \eqref{wf1} and integrate with 
respect to time getting 
\begin{gather*}
	\int_{\Omega }^{}u(t) dx+ \int_{\Gamma }^{} u_{\Gamma }(t) d\Gamma 
	= m_0 \bigl( |\Omega | + |\Gamma | \bigr) 
	\quad \mbox{for all }t \in [0,T],
\end{gather*}
where $|\Omega |:=\int_{\Omega }^{}1 dx$, $|\Gamma |:=\int_{\Gamma }^{}1 d\Gamma$ and 
\begin{gather*}
	m_0:= 
	\frac{\displaystyle \int_{\Omega }^{}u_0 dx
	+ \int_{\Gamma }^{} u_{0\Gamma } d\Gamma }{|\Omega |+|\Gamma| }
\end{gather*}
is a sort of mean value for our problem. 
Then the mean value of the variable $\mbox{\boldmath $ u$}$ is conserved 
in the sense that 
\begin{gather*}
	m \bigl( \mbox{\boldmath $ u$}(t) \bigr)=m(\mbox{\boldmath $ u$}_0) =m_0 
	\quad \mbox{for all } t \in [0,T],
\end{gather*}
where 
\begin{gather}
	m(\mbox{\boldmath $ z$}):=\frac{\displaystyle \int_{\Omega }^{}z dx
	+ \int_{\Gamma }^{} z_{\Gamma } d\Gamma }{|\Omega |+|\Gamma| }
	\quad \mbox{for all }\mbox{\boldmath $ z$} \in \mbox{\boldmath $ H$}.
\label{pier2}
\end{gather}
The duality pairing between $\mbox{\boldmath $ V$}^*$ and 
$\mbox{\boldmath $ V$}$ is denoted {\pier 
by $\langle \cdot ,\cdot \rangle _{\mbox{\scriptsize \boldmath $ V$}^*, \mbox{\scriptsize 
\boldmath $ V$}}$ and it is understood that   $\mbox{\boldmath $ H$}$  is embedded in  $
\mbox{\boldmath $ V$}^*$  in the usual way, i.e., such
that $  \langle\mbox{\boldmath $ u $},\mbox{\boldmath $ z $} \rangle_{\mbox{\scriptsize 
\boldmath $ V$}^*, \mbox{\scriptsize \boldmath $ V$}}= (\mbox{\boldmath $ u $},
\mbox{\boldmath $ z $} )_{\mbox{\scriptsize \boldmath $ H$}}$ for all  $ \mbox{\boldmath $ u 
$}  \in   \mbox{\boldmath $ H$} $ and 	$ \mbox{\boldmath $ z $} \in  \mbox{\boldmath $ V $}$.} 
Then, with the help of {\pier the later} Remark~2 (see also {\pier the comments in \cite[pp.~5674--5675]{Kub12}}), we can rewrite \eqref{wf1}~as 
\begin{gather*} 
	\bigl \langle \mbox{\boldmath $ u$}'(t),\mbox{\boldmath $ z$} 
	\bigr \rangle _{\mbox{\scriptsize \boldmath $ V$}^*, 
	\mbox{\scriptsize \boldmath $ V$}}
	+a\bigl( \mbox{\boldmath $ \mu $}(t),\mbox{\boldmath $ z$} \bigr) 
	=0 
	\quad \mbox{for all } \mbox{\boldmath $ z$} \in \mbox{\boldmath $ V$},
\end{gather*}
where $\mbox{\boldmath $ u$}'(t)$ denotes now the time derivative of the 
vectorial function and {\pier the} bilinear form 
$a(\cdot ,\cdot ):\mbox{\boldmath $ V$} \times \mbox{\boldmath $ V$} \to \mathbb{R}$ 
{\pier is defined} by 
\begin{gather*}
	a(\mbox{\boldmath $ u$},\mbox{\boldmath $ z$} ):=
	\int_{\Omega }^{} \nabla u \cdot \nabla z dx 
	+\int_{\Gamma }^{} \nabla _\Gamma u _\Gamma \cdot \nabla _\Gamma z_\Gamma d\Gamma 
	\quad \mbox{for all }\mbox{\boldmath $ u$},\mbox{\boldmath $ z$} \in \mbox{\boldmath $ V$}.
\end{gather*} 
We also {\pier introduce} the subspace $\mbox{\boldmath $ H$}_0$ of $\mbox{\boldmath $ H$}$ by 
\begin{gather*}
	\mbox{\boldmath $ H$}_0:=\left\{ \mbox{\boldmath $ z$} \in
	\mbox{\boldmath $ H$} \ : \ m(\mbox{\boldmath $ z$})=0 \right\},
\end{gather*}
and $\mbox{\boldmath $ V$}_0 :=\mbox{\boldmath $ V$} \cap \mbox{\boldmath $ H$}_0$ with 
their norms{\pier :} $| \mbox{\boldmath $ z$}|_{\mbox{\scriptsize \boldmath $ H$}_0}:=|\mbox{\boldmath $ z$}|_{\mbox{\scriptsize \boldmath $ H$}}$ for all 
$\mbox{\boldmath $ z$} \in \mbox{\boldmath $ H$}_0$ and
\begin{gather*}
	|\mbox{\boldmath $ z$}|_{\mbox{\scriptsize \boldmath $ V$}_0}:=
	{\pier \sqrt{
	a(\mbox{\boldmath $ z$},\mbox{\boldmath $ z$} )
	}}
	\quad \mbox{for all }\mbox{\boldmath $ z$} \in \mbox{\boldmath $ V$}_0.
\end{gather*}
{\pier Let us define the linear bounded} operator $\mbox{\boldmath $ F$}: \mbox{\boldmath $ V$}_0 \to \mbox{\boldmath $ V$}_0^*$ by {\pier
\begin{gather}
	\langle \mbox{\boldmath $ F$} 
	\mbox{\boldmath $ z$}, \tilde{\mbox{\boldmath $ z$}} 
	\rangle _{\mbox{\scriptsize \boldmath $ V$}_0^*, 
	\mbox{\scriptsize \boldmath $ V$}_0}
	:= a(\mbox{\boldmath $ z$},\tilde{\mbox{\boldmath $ z$}}) {\pier , 
	\quad \mbox{\boldmath $ z$}, \tilde{\mbox{\boldmath $ z$}} \in \mbox{\boldmath $ V$}_0,}
\label{pier1}
\end{gather}
as well.} Then we see that there exists $c_p>0$ such that 
\begin{gather}
	c_p |\mbox{\boldmath $ z$}|_{\mbox{\scriptsize \boldmath $ V$}}^2 
	\le \langle \mbox{\boldmath $ F$} \mbox{\boldmath $ z$}, \mbox{\boldmath $ z$} 
	\rangle _{\mbox{\scriptsize \boldmath $ V$}^*_0, 
	\mbox{\scriptsize \boldmath $ V$}_0}=|\mbox{\boldmath $ z$}|_{\mbox{\scriptsize \boldmath $ V$}_0}^2
	\quad \mbox{for all }\mbox{\boldmath $ z$} \in \mbox{\boldmath $ V$}_0 \, {\pier ;}
	\label{poin}
\end{gather}
this is checked in the Appendix. 
Therefore, thanks to the fact $|\mbox{\boldmath $ z$}|_{\mbox{\scriptsize \boldmath $ V$}_0}^2 
\le |\mbox{\boldmath $ z$}|_{\mbox{\scriptsize \boldmath $ V$}}^2$ for all 
$\mbox{\boldmath $ z$} \in \mbox{\boldmath $ V$}_0$, we see that 
$|\cdot |_{\mbox{\scriptsize \boldmath $ V$}_0}$ and $|\cdot |_{\mbox{\scriptsize \boldmath $ V$}}$ are 
equivalent norm on $\mbox{\boldmath $ V$}_0$ and then $\mbox{\boldmath $ F$}$ is the duality mapping from 
$\mbox{\boldmath $ V$}_0$ to $\mbox{\boldmath $ V$}_0^*$. 
Additionally, we can define the inner product {\pier in $\mbox{\boldmath $ V$}_0^*$}   by 
\begin{gather*}
	(\mbox{\boldmath $ z$}_1^*,\mbox{\boldmath $ z$}_2^*)_{\mbox{\scriptsize \boldmath $ V$}_0^*}
	:=\langle \mbox{\boldmath $ z$}_1^*, 
	\mbox{\boldmath $ F$} ^{-1} \mbox{\boldmath $ z$}_2^* 
	\rangle _{\mbox{\scriptsize \boldmath $ V$}^*_0, \mbox{\scriptsize \boldmath $ V$}_0}
	\quad \mbox{for all } \mbox{\boldmath $ z$}_1^*,\mbox{\boldmath $ z$}_2^* \in \mbox{\boldmath $ V$}_0^*.
\end{gather*}
Then, we obtain 
$\mbox{\boldmath $ V$}_0 
\mathop{\hookrightarrow} \mathop{\hookrightarrow}
\mbox{\boldmath $ H$}_0 
\mathop{\hookrightarrow} \mathop{\hookrightarrow}
\mbox{\boldmath $ V$}_0^*$ 
{\pier (this is also checked in the Appendix),} where 
``$\mathop{\hookrightarrow} \mathop{\hookrightarrow} $'' stands for 
the dense and compact embedding, namely 
$(\mbox{\boldmath $ V$}_0,\mbox{\boldmath $ H$}_0,\mbox{\boldmath $ V$}_0^*)$ is a standard {H}ilbert triplet.

\subsection{Definition of the solution and main theorem{\pier s}}

In order to define our solution we use that 
following {\pier additional} notation: {\pier the variable $\mbox{\boldmath $ v$}
:=\mbox{\boldmath $ u$}-m_0 \mbox{\boldmath $ 1$}$ with initial value}
$\mbox{\boldmath $ v$}_0:=\mbox{\boldmath $ u$}_0-m_0 \mbox{\boldmath $ 1$}$, 
namely, $(v,v_\Gamma )=(u-m_0,u_\Gamma -m_0)$ and 
$(v_0,v_{0\Gamma })=(u_0-m_0,u_{0\Gamma }-m_0)$; 
{\pier the datum} $\mbox{\boldmath $ f$}:=(f,f_\Gamma )$; {\pier the nonlinearity} 
$\mbox{\boldmath $ \pi $}(\mbox{\boldmath $ z$}):=(\pi (z),\pi _\Gamma (z_\Gamma ))$
{\pier for $\mbox{\boldmath $ z$}=(z,z_\Gamma) \in \mbox{\boldmath $ H$}$; 
{\pier the further space} $\mbox{\boldmath $ W$}:=H^2(\Omega ) \times H^2(\Gamma )$. }

The solution is defined as follows.

\paragraph{Definition 2.1.} 
{\it The triplet $(\mbox{\boldmath $ v$}, \mbox{\boldmath $ \mu  $}, \mbox{\boldmath $ \xi$})$ 
is called the weak solution of {\rm {\pier (P)}} if 
\begin{align*}
	\mbox{\boldmath $ v$} & \in H^1(0,T;\mbox{\boldmath $ V$}_0^*) \cap L^\infty (0,T;\mbox{\boldmath $ V$}_0) \cap L^2(0,T;\mbox{\boldmath $ W$}),
	\\
	\mbox{\boldmath $ \mu $} & \in L^2(0,T;\mbox{\boldmath $ V$}),
	\\
	\mbox{\boldmath $ \xi $} & =(\xi, \xi _\Gamma ) \in L^2(0,T;\mbox{\boldmath $ H$}), \quad 
	\xi \in \beta (v+m_0) \quad \mbox{a.e.\ in }Q, \quad 
	\xi _\Gamma \in \beta _\Gamma (v_\Gamma+m_0)\quad \mbox{a.e.\ on }\Sigma
\end{align*}
and they satisfy}
\begin{gather} 
	\bigl \langle \mbox{\boldmath $ v$}'(t),\mbox{\boldmath $ z$} 
	\bigr \rangle _{\mbox{\scriptsize \boldmath $ V$}_0^*, 
	\mbox{\scriptsize \boldmath $ V$}_0}
	+a\bigl( \mbox{\boldmath $ \mu $}(t),\mbox{\boldmath $ z$} \bigr) 
	=0 
	\quad \mbox{\it for all } \mbox{\boldmath $ z$} \in \mbox{\boldmath $ V$}_0,
	\label{d1} 
	\\
	\bigl( 
	\mbox{\boldmath $ \mu $}(t),\mbox{\boldmath $ z$}
	\bigr)_{\mbox{\scriptsize \boldmath $ H$}} 
	= a\bigl( \mbox{\boldmath $ v$}(t),\mbox{\boldmath $ z$} \bigr) 
	+ \bigl( 
	\mbox{\boldmath $ \xi $}(t) 
	+ \mbox{\boldmath $ \pi$}(\mbox{\boldmath $ v$}(t)+m_0\mbox{\boldmath $ 1$})
	-\mbox{\boldmath $ f$}(t),\mbox{\boldmath $ z$}
	\bigr)_{\mbox{\scriptsize \boldmath $ H$}} 
	\quad \mbox{\it for all }\mbox{\boldmath $ z$} \in \mbox{\boldmath $ V$},	\label{d2}
\end{gather}
{\it for a.a.\ $t\in (0,T)$, and}
\begin{gather}
	\mbox{\boldmath $ v$}(0)=\mbox{\boldmath $ v$}_0 
	\quad \mbox{\it in }\mbox{\boldmath $ H$}_0.
	\label{d3}
\end{gather} 

\paragraph{Remark 1.} Thanks to the regularity 
$\mbox{\boldmath $ v$} \in  L^2(0,T;\mbox{\boldmath $ W$})$, we see that \eqref{d2} implies that 
\begin{gather*} 
	\mu =-\Delta u + \xi +\pi (u)-f 
	\quad \mbox{a.e.\ in } Q,
	\\
	\mu _\Gamma =\partial_\nu u- \Delta _\Gamma u_\Gamma +\xi _\Gamma +\pi _\Gamma (u_\Gamma )-f_\Gamma 
	\quad \mbox{a.e.\ on } \Sigma,
\end{gather*} 
with $u=v+m_0$ and $u_\Gamma =v_\Gamma +m_0$. 
\vspace{4mm}

{\pier Next, we introduce the notion of \emph{strong solution}: we ask the reader to let us use the variable $\mbox{\boldmath $ u$}$ (instead of $\mbox{\boldmath $ v$}$) here.} 

\paragraph{Definition 2.2.} 
{\it The triplet $(\mbox{\boldmath $ u$}, \mbox{\boldmath $ \mu  $}, 
\mbox{\boldmath $ \xi$})$ 
is called the strong solution of {\rm {\pier (P)}} if 
\begin{align*}
	\mbox{\boldmath $ u$} & \in W^{1,\infty }(0,T;\mbox{\boldmath $ V$}^*) \cap H^1(0,T;\mbox{\boldmath $ V$}) \cap L^\infty (0,T;\mbox{\boldmath $ W$}),
	\\
	\mbox{\boldmath $ \mu $} & \in L^\infty (0,T;\mbox{\boldmath $ V$}) \cap L^2(0,T;\mbox{\boldmath $ W$}),
	\\
	\mbox{\boldmath $ \xi $} & \in L^\infty (0,T;\mbox{\boldmath $ H$}),
\end{align*}
and they satisfy}
\begin{gather} 
	{\partial_t u}-\Delta \mu = 0 
	\quad \mbox{\it a.e.\ in }Q,
	\label{d4}
\\
	\xi \in \beta (u),
	\quad 
	\mu = -\Delta u + \xi + \pi (u) -f 
	\quad \mbox{\it a.e.\ in }Q,
	\label{d5}
\\
	u_\Gamma =u _{|_\Gamma }, \quad \mu _\Gamma =\mu _{|_\Gamma },
	\quad 
	{\partial_t u_\Gamma }+\partial_\nu \mu -\Delta _\Gamma \mu _\Gamma =0
	\quad \mbox{\it a.e.\ on }\Sigma,
	\label{d6}
\\ 
	\xi_\Gamma  \in \beta _\Gamma (u_\Gamma ),
	\quad 
	\mu _\Gamma =\partial_\nu u - \Delta _\Gamma u_\Gamma +\xi _\Gamma + \pi _\Gamma (u_\Gamma )-f_\Gamma 
	\quad \mbox{\it a.e.\ on }\Sigma,
	\label{d7}
\\
	u(0)=u_0
	\quad 
	\mbox{\it a.e.\ in }\Omega, \quad 
	u_\Gamma  (0)=u_{0\Gamma}
	\quad \mbox{\it a.e.\ on }\Gamma.
	\label{d8}
\end{gather}

The first result states the continuous dependence on the data. 
The uniqueness of the component $\mbox{\boldmath $ v$}$ of 
the solution is also guaranteed by this theorem. 
{\pier We} assume that
\begin{enumerate}
 \item[(A1)] $\mbox{\boldmath $ f$}  \in L^2(0,T;\mbox{\boldmath $ H$}) $;
 \item[(A2)] $\mbox{\boldmath $ u$}_0 := (u_0,u_{0\Gamma }) \in \mbox{\boldmath $ V$}$;
 \item[(A3)] $\beta $, $\beta _{\Gamma }$, maximal monotone graphs in 
$\mathbb{R} \times \mathbb{R}$, are the subdifferentials 
\begin{gather*}
	\beta =\partial \widehat{\beta}, \quad \beta _{\Gamma }
	=\partial \widehat{\beta }_{\Gamma }
\end{gather*}
of some proper lower semicontinuous and convex functions 
$\widehat{\beta }$ and $\widehat{\beta }_{\Gamma }: \mathbb{R} \to [0,+\infty ]$ 
satisfying $\widehat{\beta }(0)=\widehat{\beta}_{\Gamma }(0)=0$
with some 
effective domains $D(\beta )$ and $D(\beta _\Gamma)$, respectively. 
This implies that 
$0 \in \beta (0)$ and $0 \in \beta _{\Gamma }(0)$; 
 \item[(A4)]
$\pi $, $\pi _{\Gamma }: \mathbb{R} \to \mathbb{R}$ are {L}ipschitz continuous functions 
with {L}ipschitz constants $L$ and $L_{\Gamma}$, respectively; 
\end{enumerate}

Then, we obtain the following continuous dependence on the data.

\paragraph{Theorem 2.1.} 
{\it Assume {\rm (A1)}--{\rm (A4)}. For $i=1,2$ {\pier
let $(\mbox{\boldmath $ v$}^{(i)}, \mbox{\boldmath $ \mu $}^{(i)}, 
\mbox{\boldmath $ \xi $}^{(i)})$ be a weak solution of {\rm {\pier (P)}}  
corresponding to the data} $\mbox{\boldmath $ f$}^{(i)}$ and 
$\mbox{\boldmath $ v$}_0^{(i)}$. 
Then, there 
exists a positive constant $C>0$, 
depending {\pier only on $c_p$,} $L$, $L_{\Gamma}$ and $T$, such that 
\begin{align} 
	\lefteqn{ 
	\bigl| \mbox{\boldmath $ v$}^{(1)}(t)-  \mbox{\boldmath $ v$}^{(2)}(t) 
	\bigr|_{\mbox{\scriptsize \boldmath $ V$}_0^*}^2 
	+ \int_{0}^{t}\bigl| \mbox{\boldmath $ v$}^{(1)}(s)-  \mbox{\boldmath $ v$}^{(2)}(s) 
	\bigr|_{\mbox{\scriptsize \boldmath $ V$}_0}^2
	} \nonumber \\
	& \le C \left\{ 
	\bigl| \mbox{\boldmath $ v$}^{(1)}_0-  \mbox{\boldmath $ v$}^{(2)}_0 
	\bigr|_{\mbox{\scriptsize \boldmath $ V$}_0^*}^2 
	+
	\int_{0}^{t} 
	\bigl| \mbox{\boldmath $ f$}^{(1)}(s ) -\mbox{\boldmath $ f$}^{(2)}(s ) 
	\bigr|_{\mbox{\scriptsize \boldmath $ V$}^*}^2 ds
	\right\}  \quad {\pier  \hbox{for all $t \in [0,T]$.}} \label{pier4}
\end{align} 
}
\vspace{2mm}

The second result deals with the existence of the weak solution. 
To the aim, we further assume that: 
\begin{enumerate}
 \item[(A5)] {\pier either $\mbox{\boldmath $ f$} \in W^{1,1} (0,T; \mbox{\boldmath $ H$} )  $ or} $ \mbox{\boldmath $ f$} \in L^2(0,T;\mbox{\boldmath $ V$}) $;
 \item[(A6)] $D(\beta _\Gamma ) \subseteq D(\beta )$ and 
there exist positive constants $\varrho, c_0 >0$ such that 
\begin{gather} 
	\bigl |\beta ^\circ (r) \bigr| 
	\le \varrho \bigl |\beta _{\Gamma }^\circ (r) \bigr |+c_0
	\quad 
	\mbox{for all } r \in D(\beta _\Gamma );
	\label{a4}
\end{gather} 
  \item[(A7)] $m_0\in {\rm int}D(\beta _\Gamma )$ and the {\pier compatibility conditions
$\widehat \beta (u_0) \in L^1(\Omega )$, 
$\widehat \beta _{\Gamma }(u_{0\Gamma}) \in L^1(\Gamma)$} hold.
\end{enumerate}
The minimal section $\beta ^\circ $ of $\beta $ is specified by
$\beta ^\circ (r):=\{ r^* \in \beta (r) : |r^*|=\min _{s \in \beta (r)} |s| \}$
and same definition applies to $\beta _{\Gamma }^\circ $. 
These assumptions are {\pier the same as in} \cite{CC13, CF14, CGS14}.

\paragraph{Theorem 2.2.} 
{\it Under the assumptions {\rm (A2)}--{\rm (A7)}, 
there exists a weak solution of {\pier the problem~{\rm {\pier (P)}}.}}
\vspace{2mm}

About the strong solution, we {\pier refer the reader to Subsection~4.4.}

\subsection{Abstract formulation} 

In this subsection, an abstract formulation of the problem is given. 
We can write the problem as an evolution {\pier equation including} a 
subdifferential operator{\pier : here, one can find some analogies with the approach 
followed} in \cite{CF14, CF15, KN96, Kub12}.

{\pier We define the lower semicontinuous} and convex functional 
$\varphi: \mbox{\boldmath $ H$}_0 \to [0,+\infty ]$ by 
\begin{gather*}
	\varphi (\mbox{\boldmath $ z$}) 
	:= \left\{ 
	\begin{array}{l}
	\displaystyle 
	\frac{1}{2} \int_{\Omega }^{} | \nabla z |^2 dx 
	+\frac{1}{2} \int_{\Gamma }^{} 
	 |\nabla _{\Gamma } z_{\Gamma }  |^2 d\Gamma  
	\quad 
	\mbox{if } 
	\mbox{\boldmath $ z$} \in \mbox{\boldmath $ V$}_0, \vspace{2mm}\\
	+\infty \quad \mbox{otherwise}. 
	\end{array} 
	\right. 
\end{gather*}
Then we {\pier claim} that the subdifferential $\partial \varphi $ on $\mbox{\boldmath $ H$}_0$
{\pier fulfills
$\partial \varphi (\mbox{\boldmath $ z$})=(-\Delta z,\partial_\nu z-\Delta _\Gamma z_\Gamma )$ for 
$ \mbox{\boldmath $ z$} \in D(\partial \varphi )=\mbox{\boldmath $ W$} \cap \mbox{\boldmath $ V$}_0$:
this is checked precisely in the Appendix (see Lemma~C). We also note that (cf. \eqref{pier1}--\eqref{poin})
\begin{equation} 
2 \,\varphi (\mbox{\boldmath $ z$})=  a(\mbox{\boldmath $ z$}, \mbox{\boldmath $ z$}  ) =
\langle \mbox{\boldmath $ F$} \mbox{\boldmath $ z$}, \mbox{\boldmath $ z$} 
	\rangle _{\mbox{\scriptsize \boldmath $ V$}^*_0, 
	\mbox{\scriptsize \boldmath $ V$}_0}=|\mbox{\boldmath $ z$}|_{\mbox{\scriptsize \boldmath $ V$}_0}^2
	\quad \mbox{for all }\mbox{\boldmath $ z$} \in \mbox{\boldmath $ V$}_0,
	\label{pier14}
\end{equation}
by collecting in this formula part of our notation.} 

{\pier At this point, we emphasize that our problem is equivalent to the following {C}auchy problem for a suitable evolution equation}.
\begin{gather} 
	\mbox{\boldmath $ v$}'(t) + 
	{\revis 
	\mbox{\boldmath $ F$} \bigl (\mbox{\boldmath $ P$} \mbox{\boldmath $ \mu $} (t) \bigr) }
	= {\revis \mbox{\boldmath $ 0$}} \quad \mbox{in }\mbox{\boldmath $ V$}_0^*,
	\ \mbox{for a.a.\ }t\in (0,T),
	\label{e1}
	\\
	\mbox{\boldmath $ \mu $}(t)
	=
	\partial \varphi {\revis \bigl (} \mbox{\boldmath $ v$}(t) {\revis \bigr )} 
	+ \mbox{\boldmath $ \xi $}(t)
	+ \mbox{\boldmath $ \pi $} \bigl (\mbox{\boldmath $ u$}(t) \bigr ) 
	-\mbox{\boldmath $ f$}(t)
	\quad \mbox{in } \mbox{\boldmath $ H$},
	\ \mbox{for a.a.\ } t \in (0,T),
	\label{e2}
	\\
	\mbox{\boldmath $ \xi $}(t) \in \mbox{\boldmath $ \beta $}
	\bigl ( \mbox{\boldmath $ u$}(t) \bigr ) 
	\quad \mbox{in } \mbox{\boldmath $ H$},
	\ \mbox{for a.a.\ } t \in (0,T),
	\label{e3}
	\\
	\mbox{\boldmath $ v$}(0)=\mbox{\boldmath $ v$}_0 
	\quad \mbox{in } \mbox{\boldmath $ H$}_0,
	\label{ic}
\end{gather} 
where {\pier $\mbox{\boldmath $ u$}=\mbox{\boldmath $ v$}+m_0\mbox{\boldmath $ 1$}=(v+m_0,v_\Gamma +m_0)$ and $\mbox{\boldmath $ \beta $}(\mbox{\boldmath $ z$})
:=(\beta (z),\beta _\Gamma(z_\Gamma ))$, 
\begin{gather}
	\mbox{\boldmath $ P$}\mbox{\boldmath $ z$}:= 
	\mbox{\boldmath $ z$}-m(\mbox{\boldmath $ z$}) \mbox{\boldmath $ 1$}
	= \bigl (z - m(\mbox{\boldmath $ z$}), z_\Gamma  - m(\mbox{\boldmath $ z$}) \bigr ) ,
	\label{pier10}
\end{gather}
with $m(\mbox{\boldmath $ z$}) $ defined by \eqref{pier2}, 
for all $\mbox{\boldmath $ z$}  = (z , z_\Gamma ) \in \mbox{\boldmath $ H$}.$
Note that {\pier the projection operator}
$\mbox{\boldmath $ P$}$ acts as linear bounded operator both from 
$\mbox{\boldmath $ V$} $ to $\mbox{\boldmath $ V$}_0$ and
from  $\mbox{\boldmath $ H$} $ to $\mbox{\boldmath $ H$}_0$.
Moreover,  it is easy to} see that 
\begin{align}
	(\mbox{\boldmath $ z$}^*, \mbox{\boldmath $ P$}\tilde{\mbox{\boldmath $ z$}} 
	)_{\mbox{\scriptsize \boldmath $ H$}_0} 
	& = \left\{ \int_{\Omega }^{}z^* \tilde{z}dx 
	 + \int_{\Gamma }^{}z_\Gamma^* \tilde{z}_\Gamma \right\} 
	 - m(\tilde{\mbox{\boldmath $ z$}}) \left\{ \int_{\Omega }^{}z^* dx 
	 + \int_{\Gamma }^{}z_\Gamma^*  \right\}
	\nonumber \\
	& = (\mbox{\boldmath $ z$}^*,\tilde{\mbox{\boldmath $ z$}})_{\mbox{\scriptsize \boldmath $ H$}} 
	\quad \mbox{for all }\mbox{\boldmath $ z$}^* \in \mbox{\boldmath $ H$}_0 
	\ \mbox{and } \tilde{\mbox{\boldmath $ z$}} \in \mbox{\boldmath $ H$}.
	\label{proj}
\end{align}

\paragraph{Remark 2.} Note that in Definition{\pier s}~2.1 and 2.2, 
$\mbox{\boldmath $ v$}' \in L^2(0,T;\mbox{\boldmath $ V$}_0^*)$ 
can be easily extended to 
$L^2(0,T;\mbox{\boldmath $ V$}^*)$ by setting 
$\langle \mbox{\boldmath $ v$}'(t),\mbox{\boldmath $ 1$} \rangle_{\mbox{\scriptsize \boldmath $ V$}^*,\mbox{\scriptsize \boldmath $ V$}} :=0$, namely 
\begin{gather*} 
	\bigl\langle \mbox{\boldmath $ v$}'(t),\mbox{\boldmath $ z$} 
	\bigr\rangle _{\mbox{\scriptsize \boldmath $ V$}^*,\mbox{\scriptsize \boldmath $ V$}}
	:= \bigl\langle 
	\mbox{\boldmath $ v$}'(t),\mbox{\boldmath $ P$}\mbox{\boldmath $ z$}
	\bigr\rangle _{\mbox{\scriptsize \boldmath $ V$}_0^*,\mbox{\scriptsize \boldmath $ V$}_0}
	\quad \mbox{for all } \mbox{\boldmath $ z$} \in \mbox{\boldmath $ V$},
\end{gather*} 
because we know the following orthogonal decomposition 
$\mbox{\boldmath $ V$}=\mbox{\boldmath $ V$}_0 \oplus \mbox{\boldmath $ R$}$ with 
$\mbox{\boldmath $ R$}=\{ r \mbox{\boldmath $ 1$} : r \in \mathbb{R}\}$. 
Therefore, {\pier we have
\begin{align*}
\bigl \langle \mbox{\boldmath $ v$}'(t), \mbox{\boldmath $ 1$}
\bigr \rangle _{\mbox{\scriptsize \boldmath $ V$}^*,\mbox{\scriptsize \boldmath $ V$}} 
 = 
\bigl \langle \mbox{\boldmath $ v$}'(t),\mbox{\boldmath $ P$} \mbox{\boldmath $ 1$}
\bigr \rangle _{\mbox{\scriptsize \boldmath $ V$}^*_0,\mbox{\scriptsize \boldmath $ V$}_0}
 = \bigl \langle \mbox{\boldmath $ v$}'(t), \mbox{\boldmath $ 1$}- 
m(\mbox{\boldmath $ 1$})\mbox{\boldmath $ 1$}
\bigr \rangle _{\mbox{\scriptsize \boldmath $ V$}^*_0,\mbox{\scriptsize \boldmath $ V$}_0}
\end{align*}
and it is clear that} $\mbox{\boldmath $ 1$}-m(\mbox{\boldmath $ 1$})=0$ 
in $\mbox{\boldmath $ V$}_0$. 
Thus, we can identify {\pier the dual space}
$\mbox{\boldmath $ V$}_0^*$ by $\{ \mbox{\boldmath $ z$}^* \in \mbox{\boldmath $ V$}^* : 
\langle \mbox{\boldmath $ z$}^*, \mbox{\boldmath $ 1$} 
\rangle _{\mbox{\scriptsize \boldmath $ V$}^*,\mbox{\scriptsize \boldmath $ V$}}=0\}$ and  
we know that \eqref{d1} holds for all $\mbox{\boldmath $ z$} \in \mbox{\boldmath $ V$}$, that is, 
its variational formulation {\pier implies} the structure of volume conservation 
(see, {\pier e.g.,}~\cite[p.~5674]{Kub12}). 
\vspace{2mm}

Now, {\pier it is straightforward to check that \eqref{e1} and \eqref{e2} yield the single abstract equation}
\begin{align} 
	\mbox{\boldmath $ F$}^{-1} 
	\bigl( \mbox{\boldmath $ v$}'(t) \bigr) 
	+  \partial \varphi {\revis \bigl( }\mbox{\boldmath $ v$}(t) {\revis \bigr)} 
	= \mbox{\boldmath $ P$} \bigl( - \mbox{\boldmath $ \xi $}(t)
	- \mbox{\boldmath $ \pi $} \bigl ( \mbox{\boldmath $ v$}(t) 
	+ m_0 \mbox{\boldmath $ 1$} \bigr ) 
	+ \mbox{\boldmath $ f$}(t) \bigr) \quad 
	\nonumber
	\\ 
	\mbox{in } \mbox{\boldmath $ H$}_0, \ \mbox{for a.a.\ } t \in (0,T),
	\label{e4}
\end{align}
{\pier where 
\begin{gather}
\mbox{\boldmath $ \xi $}(t) \in 
	- \mbox{\boldmath $ \beta $} 
	\bigl (\mbox{\boldmath $ v$}(t) 
	+ m_0 \mbox{\boldmath $ 1$}\bigr )
	\quad \mbox{in } \mbox{\boldmath $ H$}, \ \mbox{for a.a.\ } t \in (0,T).
\label{pier3}
\end{gather}
Based on the abstract theory {\pier developed in}~\cite{CV90}, 
we expect that \eqref{e4}--\eqref{pier3} can be solved by the approach for} doubly nonlinear evolution inclusion{\pier s,} which will be discussed in Section 4. 

\section{Continuous dependence}
\setcounter{equation}{0}

In this section, we prove Theorem~2.1. 

\paragraph{Proof of Theorem 2.1.} For $i=1,2$ let 
$(\mbox{\boldmath $ v$}^{(i)}, \mbox{\boldmath $ \mu $}^{(i)}, 
\mbox{\boldmath $ \xi $}^{(i)})$
be a weak solution of {\rm {\pier (P)}} corresponding to the data 
($\mbox{\boldmath $ f$}^{(i)}$, $\mbox{\boldmath $ v$}_{0}^{(i)}$). 
We consider the difference between {\pier equations \eqref{e4} written, 
at the time $s \in (0,T)$, for 
$\mbox{\boldmath $ v$}^{(1)}(s)=(v^{(1)}(s),v_\Gamma ^{(1)}(s))$ 
and 
$\mbox{\boldmath $ v $}^{(2)}(s)=(v^{(2)}(s),v_\Gamma ^{(2)}(s))$,
respectively.}
Then, we take 
the inner product with
$\mbox{\boldmath $ v$}^{(1)}(s)-\mbox{\boldmath $ v$}^{(2)}(s)$ in 
$\mbox{\boldmath $ H$}_0$. 
Using the property 
\eqref{proj}  {\pier as well as \eqref{pier3} and} the monotonicity of 
$\beta $, $\beta _\Gamma $
we obtain 
\begin{align} 
	\lefteqn{ 
	\frac{1}{2} \frac{d}{ds} 
	\bigl| \mbox{\boldmath $ v$}^{(1)}(s)-  \mbox{\boldmath $ v$}^{(2)}(s) 
	\bigr|_{\mbox{\scriptsize \boldmath $ V$}_0^*}^2 
	+ \bigl| \mbox{\boldmath $ v$}^{(1)}(s)-  \mbox{\boldmath $ v$}^{(2)}(s) 
	\bigr|_{\mbox{\scriptsize \boldmath $ V$}_0}^2
	} \nonumber \\
	& \le {} {\revis \bigl \langle} 
	\mbox{\boldmath $ f$}^{(1)}(s) -\mbox{\boldmath $ f$}^{(2)}(s),
	\mbox{\boldmath $ v$}^{(1)}(s) -  \mbox{\boldmath $ v$}^{(2)}(s)  
	{\revis \bigr \rangle}_{\mbox{\scriptsize \boldmath $ V$}^*, \mbox{\scriptsize \boldmath $ V$}}
	\nonumber \\
	& \quad {}
	- \Bigl( \pi \bigl (v^{(1)}(s) +m_0\bigr) -  \pi \bigl (v^{(2)}(s) +m_0\bigr) , 
	v^{(1)}(s)-  v^{(2)}(s) \Bigr)_{\! H}
	\nonumber \\
	& \quad {} - \Bigl( \pi_\Gamma  \bigl (v^{(1)}_\Gamma (s) +m_0 \bigr) 
	-  \pi_\Gamma  \bigl (v^{(2)}_\Gamma (s) +m_0\bigr) , 
	v^{(1)}_\Gamma (s)-  v^{(2)}_\Gamma (s) \Bigr)_{\! H_\Gamma }
	\label{1st}
\end{align} 
for a.a.\ $s \in (0,T)$. 
By virtue of the compact {\pier embedding of $\mbox{\boldmath $ V$}_0$ into
$\mbox{\boldmath $ H$}_0$, for each $\delta >0$ there exists a positive 
constant $c_\delta$ (depending on $\delta $)} such that 
\begin{gather} 
	|\mbox{\boldmath $ z$}|_{\mbox{\scriptsize \boldmath $ H$}_0}
	\le \delta |\mbox{\boldmath $ z$}|_{\mbox{\scriptsize \boldmath $ V$}_0} 
	+ c_\delta |\mbox{\boldmath $ z$}|_{\mbox{\scriptsize \boldmath $ V$}_0^*}
	\quad {\rm for~all}~ \mbox{\boldmath $ z$} \in \mbox{\boldmath $ V$}_0, 
	\label{compin}
\end{gather} 
(see, e.g., \cite[{\pier Sect.~8, Lemma~8}]{Sim87}). 
Then, {\pier in view of} \eqref{poin} and the {L}ipschitz continuities of 
$\pi$ and $\pi _\Gamma$ we {\pier infer that}
\begin{align*} 
	\lefteqn{ 
	\frac{1}{2}\frac{d}{ds} 
	\bigl| \mbox{\boldmath $ v$}^{(1)}(s)-  \mbox{\boldmath $ v$}^{(2)}(s) 
	\bigr|_{\mbox{\scriptsize \boldmath $ V$}_0^*}^2 
	+ c_p \bigl| \mbox{\boldmath $ v$}^{(1)}(s)-  \mbox{\boldmath $ v$}^{(2)}(s) 
	\bigr|_{\mbox{\scriptsize \boldmath $ V$}}^2
	} \nonumber \\
	& \le \frac{1}{c_p}
	\bigl| \mbox{\boldmath $ f$}^{(1)}(s ) -\mbox{\boldmath $ f$}^{(2)}(s ) 
	\bigr|_{\mbox{\scriptsize \boldmath $ V$}^*}^2 +
	\frac{c_p}{4} \bigl| \mbox{\boldmath $ v$}^{(1)}(s)-  \mbox{\boldmath $ v$}^{(2)}(s) 
	\bigr|_{\mbox{\scriptsize \boldmath $ V$}}^2 
	\nonumber \\
	& \quad {} 
	+ (L+L_\Gamma ) \left\{ 
	2\delta 
	\bigl| 
	\mbox{\boldmath $ v$}^{(1)}(s)-  \mbox{\boldmath $ v$}^{(2)}(s)
	\bigr|_{\mbox{\scriptsize \boldmath $ V$}_0}^2
	+ 2c_\delta 
	\bigl| 
	\mbox{\boldmath $ v$}^{(1)}(s)-  \mbox{\boldmath $ v$}^{(2)}(s)
	\bigr|_{\mbox{\scriptsize \boldmath $ V$}_0^*}^2
	\right\}, 
\end{align*} 
for a.a.\  $s \in (0,T)$, 
Therefore, taking $\delta =c_p/(8(L+L_\Gamma))$
and applying the {G}ronwall lemma{\pier,  we find a constant $C>0$, 
with the dependencies specified in the statement, such that \eqref{pier4} holds.} \hfill $\Box$

\section{Existence}
\setcounter{equation}{0}

This section is devoted to the proof of Theorem 2.2. 
We make use of the Yosida approximation for maximal monotone operators 
$\beta $, $\beta _\Gamma $ 
and of well-known results of this theory ({\pier see}~\cite{Bar10, Bre73, Ken07}). 
For each $\varepsilon \in (0,1]$, we define 
$\beta _\varepsilon, \beta _{\Gamma,\varepsilon}:\mathbb{R} \to \mathbb{R}$, 
along with the associated resolvent operators 
$J_\varepsilon, J_{\Gamma,\varepsilon}:\mathbb{R} \to \mathbb{R}$
by 
\begin{gather*}
	\beta _\varepsilon (r)
	:= \frac{1}{\varepsilon } \bigl( r-J_\varepsilon (r) \bigr)
	:=\frac{1}{\varepsilon }\bigl( r-(I+\varepsilon \beta )^{-1} (r)\bigr),
	\\
	\beta _{\Gamma, \varepsilon} (r)
	:= \frac{1}{\varepsilon \varrho} \bigl( r-J_{\Gamma,\varepsilon  }(r) \bigr )
	:=\frac{1}{\varepsilon \varrho}\bigl( r- (I+\varepsilon \varrho \beta _\Gamma )^{-1} (r) \bigr)
\end{gather*}
{\pier for all $ r \in \mathbb{R}$,}
where $\varrho>0$ is {\pier the} same constant as in the assumption \eqref{a4}. 
Note that the two definitions are not symmetric since in the second it is 
$\varepsilon \varrho$ and not directly $\varepsilon $ to be used as approximation parameter. 
Now, we easily have 
$\beta _\varepsilon(0)=\beta _{\Gamma, \varepsilon}(0)=0$. 
Moreover, the related {M}oreau--{Y}osida regularizations $\widehat{\beta }_\varepsilon, 
\widehat{\beta }_{\Gamma,\varepsilon}$
of $\widehat{\beta }, \widehat{\beta}_{\Gamma }:\mathbb{R} \to \mathbb{R}$ fulfill
\begin{gather*}
	\widehat{\beta }_{\varepsilon }(r)
	:=\inf_{s \in \mathbb{R}}\left\{ \frac{1}{2\varepsilon } |r-s|^2
	+\widehat{\beta }(s) \right\} 
	= 
	\frac{1}{2\varepsilon } 
	\bigl| r-J_\varepsilon (r) \bigr|^2+\widehat{\beta }\bigl (J_\varepsilon (r) \bigr )
	= \int_{0}^{r} \beta _\varepsilon (s)ds,
	\\
	\widehat{\beta }_{\Gamma, \varepsilon }(r)
	:=\inf_{s \in \mathbb{R}}\left\{ \frac{1}{2\varepsilon \varrho } |r-s|^2+\widehat{\beta }_\Gamma (s) \right\} 
	= \int_{0}^{r} \beta _{\Gamma,\varepsilon} (s)ds
	\quad \mbox{for all } r\in \mathbb{R}.
\end{gather*}
It is well known that $\beta_\varepsilon$ is {L}ipschitz continuous with {L}ipschitz constant 
$1/\varepsilon $ 
and $\beta_{\Gamma, \varepsilon}$ is also {L}ipschitz continuous with constant 
$1/(\varepsilon \varrho)$. In addition, we have the standard properties{\pier
\begin{align}
	&\bigl |\beta _\varepsilon (r) \bigr | \le \bigl |\beta ^\circ (r) \bigr |, \quad 
	\bigl |\beta _{\Gamma ,\varepsilon }(r) \bigr | \le \bigl |\beta _{\Gamma }^\circ (r) \bigr | 
	\quad \mbox{and }
	\nonumber \\
	&\quad 0 \le \widehat{\beta }_\varepsilon (r) \le \widehat{\beta }(r), \quad 
	0 \le \widehat{\beta }_{\Gamma, \varepsilon} (r) \le \widehat{\beta }_{\Gamma }(r)
	\quad \mbox{for all } r \in \mathbb{R}.
	\label{prim}
\end{align}
}%
Here, {\pier thanks to \cite[Lemma~4.4]{CC13}, we have that}
\begin{gather}
	\bigl |\beta_\varepsilon (r)\bigr | 
	\le \varrho \bigl |\beta _{\Gamma,\varepsilon} (r)\bigr |+c_0
	\quad 
	\mbox{for all } r \in \mathbb{R},
	\label{a4e}
\end{gather} 
{\pier with the same constants $ \varrho $ and $c_0 $ as in \eqref{a4}.}

\subsection{Approximation of the problem}

In this subsection, 
we consider an approximation for \eqref{e4} stated as the following 
{C}auchy problem: {\pier for} each $\varepsilon \in (0,1]$ find 
$\mbox{\boldmath $ v$}_\varepsilon :=(v_\varepsilon, v_{\Gamma, \varepsilon})$ 
satisfying
\begin{align} 
	&\varepsilon \mbox{\boldmath $ v$}_\varepsilon '(t)
	+ \mbox{\boldmath $ F$}^{-1}\mbox{\boldmath $ v$}_\varepsilon '(t)
	+\partial \varphi 
	\bigl( \mbox{\boldmath $ v$}_\varepsilon (t) \bigr) 
	\nonumber 
	\\
	&= \mbox{\boldmath $ P$} \bigl( 
	- \mbox{\boldmath $ \beta $}_\varepsilon 
	\bigl ( {\pier \mbox{\boldmath $ v$}_\varepsilon (t)    
	 +m_0 \mbox{\boldmath $ 1$} } \bigr )
	- \mbox{\boldmath $ \pi $} \bigl ( 
	{\pier \mbox{\boldmath $ v$}_\varepsilon (t)    
	 +m_0 \mbox{\boldmath $ 1$} } \bigr )
	+ \mbox{\boldmath $ f$}(t) \bigr) 
	\quad \mbox{in } \mbox{\boldmath $ H$}_0, 
	\ \mbox{for a.a.\ } t \in (0,T),
	\label{ape}
	\\
	&\hspace{6cm}\mbox{\boldmath $ v$}_\varepsilon (0)=\mbox{\boldmath $ v$}_0 
	\quad \mbox{in } \mbox{\boldmath $ H$}_0.
	\label{ice}
\end{align} 
The structure of this approximate problem {\pier fits into the framework 
of the general problem treated in~\cite{CV90}.}  Namely, thanks to the abstract theory 
of doubly nonlinear evolution inclusion{\pier s}, we can solve the {\pier {C}auchy problem} \eqref{ape}--\eqref{ice}. 

\paragraph{Proposition 4.1.} 
{\it 
For each $\varepsilon \in (0,1]$, there 
exists a unique 
\begin{gather*}
	\mbox{\boldmath $ v$}_\varepsilon \in H^1(0,T;\mbox{\boldmath $ H$}_0) 
	\cap 
	C {\revis \bigl (}[0,T];\mbox{\boldmath $ V$}_0 {\revis \bigr)} \cap L^2(0,T;\mbox{\boldmath $ W$}),
\end{gather*} 
such that $\mbox{\boldmath $ v$}_\varepsilon $ satisfies \eqref{ape}--\eqref{ice}. }

\paragraph{Proof.} {\pier As the argumentation follows the analogous proof performed in~\cite{CF15},
we only sketch it.} 
We claim that for a given 
$\bar{\mbox{\boldmath $ v$}}:=(\bar{v},\bar{v}_\Gamma) \in C([0,T];\mbox{\boldmath $ H$}_0)$ 
there exists a unique 
$${\pier \mbox{\boldmath $ v$} \in  H^1(0,T;\mbox{\boldmath $ H$}_0) 
	\cap 
	L^\infty (0,T;\mbox{\boldmath $ V$}_0) 
	\subset C{\revis \bigl (}[0,T];\mbox{\boldmath $ H$}_0 {\revis \bigr )}}$$
such that{\pier
\begin{align*}
	&\bigl (\varepsilon  \mbox{\boldmath $ I$} 
	+\mbox{\boldmath $ F$}^{-1} \bigr) \mbox{\boldmath $ v$} '(t)
	+\partial \varphi
	\bigl( \mbox{\boldmath $ v$} (t) \bigr) 
	\\
	&\ni \mbox{\boldmath $ P$} \bigl( 
	- \mbox{\boldmath $ \beta $}_\varepsilon 
	\bigl ( \bar{\mbox{\boldmath $ v$}} (t) 
	+m_0 \mbox{\boldmath $ 1$} \bigr )
	- \mbox{\boldmath $ \pi $} \bigl ( \bar{\mbox{\boldmath $ v$}}(t)  
	+m_0 \mbox{\boldmath $ 1$} \bigr )
	+ \mbox{\boldmath $ f$}(t) \bigr) 
	\quad \mbox{in } \mbox{\boldmath $ H$}_0,
	\ \mbox{for a.a.\ } t \in (0,T),
	\\
	&\hspace{5.5cm}  \mbox{\boldmath $ v$}(0)=\mbox{\boldmath $ v$}_0 
	\quad \mbox{in } \mbox{\boldmath $ H$}_0.
\end{align*}
}
Indeed, it suffices to apply the abstract theory of doubly nonlinear evolution inclusion{\pier s}
(see {\pier in particular} \cite[Thm.~2.1]{CV90}). 
We point out that, thanks to {\pier the presence  of $\varepsilon \mbox{\boldmath $ I$}$, the operator $\varepsilon \mbox{\boldmath $ I$}  +\mbox{\boldmath $ F$}^{-1}$ is coercive in $\mbox{\boldmath $ H$}_0 $, which is an important assumption of Theorem~2.1 in} \cite{CV90}. 
Then, we construct the map 
$\Psi : \bar{\mbox{\boldmath $ v$}} \mapsto \mbox{\boldmath $ v$}$ 
from $C([0,T];\mbox{\boldmath $ H$}_0)$ into itself. Next, 
for {\pier a} given 
$\bar{\mbox{\boldmath $ v$}}^{(i)} \in C([0,T];\mbox{\boldmath $ H$}_0)$,  {\pier we}
put $\mbox{\boldmath $ v$}^{(i)}:=\Psi \bar{\mbox{\boldmath $ v$}}^{(i)}$, $i=1,2$. 
{\pier Using} the monotonicity of $\partial \varphi$, it is not difficult to deduce the
estimate
\begin{gather*}
	\bigl| \mbox{\boldmath $ v$}^{(1)}(t) 
	- \mbox{\boldmath $ v$}^{(2)}(t) \bigr|_{\mbox{\boldmath \scriptsize $ H$}_0}^2
	\le c_{\varepsilon }
	\int_{0}^{t} 
	\bigl| \bar{\mbox{\boldmath $ v$}}^{(1)}(s ) 
	-\bar{\mbox{\boldmath $ v$}}^{(2)}(s ) \bigr|_{\mbox{\boldmath \scriptsize $ H$}_0}^2 ds 
	\quad \mbox{for all } t \in [0,T],
\end{gather*} 
where $c_\varepsilon $ is a constant depending on $L$, $L_\Gamma$ and $\varepsilon $. 
Therefore, we can prove that there exist{\pier s} {\pier a} suitable 
$k \in \mathbb{N}$ such that $\Psi ^k$ is a 
contraction mapping in $C([0,T];\mbox{\boldmath $ H$}_0)$. Hence,
being $\varepsilon >0$ there exists a unique fixed point for $\Psi$ which yields
the unique solution $\mbox{\boldmath $ v$}_\varepsilon $ of the 
problem \eqref{ape}--\eqref{ice}. 
Finally, thanks to the fact {\pier that
$\partial \varphi (\mbox{\boldmath $ v$}_\varepsilon ) \in L^2(0,T;\mbox{\boldmath $ H$}_0)$, we easily reach the conclusion, by simply observing that $H^1(0,T;\mbox{\boldmath $ H$}_0) 
	\cap L^2(0,T;\mbox{\boldmath $ W$})$ is contained in $
	C ([0,T];\mbox{\boldmath $ V$}_0) $.}
\hfill $\Box$ \\
\vspace{-3mm}

Now, for each $\varepsilon \in (0,1]$ {\pier we set} 
\begin{gather} 
	\mbox{\boldmath $ \mu $}_\varepsilon (t):=
	\varepsilon \mbox{\boldmath $ v$}_\varepsilon '(t)
	+ \partial \varphi \bigl (\mbox{\boldmath $ v$}_\varepsilon (t) \bigr )
	+ \mbox{\boldmath $ \beta $}_\varepsilon 
	\bigl (\mbox{\boldmath $ u$}_\varepsilon (t) \bigr ) 
	+ \mbox{\boldmath $ \pi $}
	\bigl (\mbox{\boldmath $ u$}_\varepsilon (t) \bigr ) 
	-\mbox{\boldmath $ f$}(t)
	\quad \mbox{for a.a.\ } t \in (0,T),
	\label{mue}
\end{gather}
where {\pier $\mbox{\boldmath $ u$}_\varepsilon :=
\mbox{\boldmath $ v$}_\varepsilon    
	 +m_0 \mbox{\boldmath $ 1$} = (v_\varepsilon +m_0,v_{\Gamma,\varepsilon  }+m_0)$.}
We expect {\pier \eqref{mue} to give} the approximate sequence for the chemical potential. 
Then, we can rewrite the evolution equation \eqref{ape}~as 
\begin{gather*}
	\mbox{\boldmath $ F$}^{-1}\bigl( \mbox{\boldmath $ v$}_\varepsilon '(t) \bigr) 
	+\mbox{\boldmath $ \mu $}_\varepsilon (t)- \omega_\varepsilon (t)  {\pier \/
	 \mbox{\boldmath $ 1$}}
	= {\revis \mbox{\boldmath $ 0$}} \quad \mbox{in }\mbox{\boldmath $ V$}, 
	\ \mbox{for a.a.\ }t \in (0,T),
\end{gather*}
where 
\begin{gather}
	\omega _\varepsilon (t)
	:= m \bigl(  \mbox{\boldmath $ \beta $}_\varepsilon 
	\bigl (\mbox{\boldmath $ u$}_\varepsilon (t) \bigr) 
	+ \mbox{\boldmath $ \pi $}
	\bigl (\mbox{\boldmath $ u$}_\varepsilon (t) \bigr) 
	-\mbox{\boldmath $ f$}(t) \bigr) ,
	{\pier \quad t \in (0,T).}
	\label{ome}
\end{gather}
{\pier Hence, we realize} that $ \mbox{\boldmath $ P$} \mbox{\boldmath $ \mu $}_\varepsilon 
= \mbox{\boldmath $ \mu $}_\varepsilon - \omega_\varepsilon \mbox{\boldmath $ 1$} \in L^2(0,T;\mbox{\boldmath $ V$}_0)$ 
and $\omega _\varepsilon \in L^2(0,T)$. {\pier From these regularities it follows that}
$\mbox{\boldmath $ \mu $}_\varepsilon \in L^2(0,T;\mbox{\boldmath $ V$})$ and {\pier the 
pair
$(\mbox{\boldmath $ v$}_\varepsilon , \mbox{\boldmath $ \mu $}_\varepsilon)$ satisfies}
\begin{gather} 
	\mbox{\boldmath $ v$}_\varepsilon '(t) + \mbox{\boldmath $ F$} 
	\bigl( \mbox{\boldmath $ P$} \mbox{\boldmath $ \mu $}_\varepsilon (t)\bigr) 
	= {\revis \mbox{\boldmath $ 0$}}
	\quad \mbox{in }\mbox{\boldmath $ V$}_0^*, 
	\ \mbox{for a.a.\ } t\in (0,T).
	\label{eonvs}
\end{gather} 
Moreover, we see that $(v_\varepsilon ,v_{\Gamma,\varepsilon}, \mu _\varepsilon , 
\mu _{\Gamma,\varepsilon  })$ {\pier fulfill} the following weak formulations: 
\begin{multline}
	\lefteqn{ 
	\int_{\Omega }^{}  {\partial_t v_\varepsilon }(t) z dx 
	+ \int_{\Gamma }^{} {\partial_t v_{\Gamma,\varepsilon } }(t) z_\Gamma d \Gamma 
	+ \int_{\Omega }^{} \nabla \mu _\varepsilon (t) \cdot \nabla z dx 
	} \\
	{\pier {}+ \int_{\Gamma }^{} \nabla_\Gamma  \mu _{\Gamma,\varepsilon }(t) \cdot \nabla_\Gamma z_\Gamma d\Gamma = 0 }
	\quad \mbox{for all } 
	\mbox{\boldmath $ z$}{\pier {}=(z, z_\Gamma)} \in \mbox{\boldmath $ V$}_0,
	\label{wf1e}
\end{multline} 
\begin{align} 
	&\int_{\Omega }^{} \mu_\varepsilon (t) z dx 
	+ \int_{\Gamma }^{} \mu _{\Gamma,\varepsilon }(t) z_\Gamma d \Gamma 
 \nonumber \\
	& = 
	\varepsilon \int_{\Omega }^{}  {\partial_t v_\varepsilon }(t) z dx 
	+ \varepsilon \int_{\Gamma }^{} {\partial_t v_{\Gamma,\varepsilon } } (t) z_\Gamma d \Gamma 
	+ \int_{\Omega }^{} \nabla v_\varepsilon (t)\cdot \nabla z dx + 
	\int_{\Gamma }^{} \nabla _\Gamma v_{\Gamma,\varepsilon }(t) \cdot \nabla _\Gamma z_\Gamma d\Gamma 
	\nonumber \\
	  &\quad {} + \int_{\Omega }^{} \bigl( 
	\beta_\varepsilon  \bigl( v_\varepsilon (t)+m_0 \bigr) 
	+\pi \bigl( v_\varepsilon (t)+m_0 \bigr)-f(t) \bigr)z dx 
	\nonumber \\
	&\quad {} + \int_{\Gamma }^{} \bigl( 
	\beta_{\Gamma, \varepsilon} 
	\bigl( v_{\Gamma,\varepsilon } (t)+m_0 \bigr) 
	+\pi_\Gamma  \bigl( v_{\Gamma,\varepsilon } (t)+m_0 \bigr)
	-f_\Gamma(t) \bigr )z_\Gamma  d\Gamma 
	\nonumber
	\\
	& \hspace{10cm}\mbox{for all } \mbox{\boldmath $ z$} {\pier {}=(z, z_\Gamma)}\in \mbox{\boldmath $ V$},
	\label{wf2e}
\end{align}
{\pier for a.a.\ $t\in (0,T)$}. 
In particular, from {\pier \eqref{wf2e} we} deduce the equations
\begin{gather} 
	\mu_\varepsilon  = \varepsilon {\partial_t v_\varepsilon }
	-\Delta v_\varepsilon  + \beta _\varepsilon (v_\varepsilon +m_0) + 
	\pi (v_\varepsilon +m_0) -f 
	\quad \mbox{a.e.\ in }Q,
	\label{d5e}
\\
	\mu _{\Gamma,\varepsilon } = 
	\varepsilon {\partial_t v_{\Gamma,\varepsilon  }}
	+ \partial_\nu v_\varepsilon  
	- \Delta _\Gamma v_{\Gamma,\varepsilon } 
	+ \beta _{\Gamma,\varepsilon }(v_{\Gamma,\varepsilon  }+m_0) 
	+ \pi _\Gamma (v_{\Gamma,\varepsilon } +m_0)
	- f_{\Gamma} 
	\quad \mbox{a.e.\ on }\Sigma.
	\label{d7e}
\end{gather} 
{\pier Arguing on \eqref{wf1e}}, since ${\pier \partial_t v_\varepsilon}\in L^2(0,T;H)$
and ${\pier \partial_t v_{\Gamma,\varepsilon}} \in L^2(0,T;H_\Gamma )$, 
we can recover that 
$\Delta \mu _\varepsilon \in L^2(0,T;H)$. 
On the other hand, we already know that 
$\mu _\varepsilon \in L^2 (0,T;V)$ and $\mu_{\Gamma, \varepsilon } \in L^2(0,T;V_\Gamma )$. 
Then, we infer that (see, e.g., \cite[Thm.~3.2, p.~1.79]{BG87})
\begin{gather*}
	\mu _\varepsilon \in L^2 \bigl (0,T;H^{3/2} (\Omega ) \bigr )
\end{gather*}
and consequently, by a trace theorem \cite[Thm.~2.27, p.~1.64]{BG87}, 
$\partial_\nu \mu_\varepsilon \in L^2(0,T;H_\Gamma )$.
Therefore, we also obtain that 
$\Delta _\Gamma \mu _{\Gamma ,\varepsilon } \in L^2(0,T;H_\Gamma )$  
{\pier so that we can write} the equations 
\begin{gather} 
	{\partial_t v_\varepsilon }-\Delta \mu_\varepsilon  = 0 
	\quad \mbox{a.e.\ in }Q,
	\label{d4e}
\\
	{\partial_t v_{\Gamma ,\varepsilon } }
	+ \partial_\nu \mu_\varepsilon  
	- \Delta _\Gamma \mu _{\Gamma ,\varepsilon } =0
	\quad \mbox{a.e.\ on }\Sigma.
	\label{d6e}
\end{gather} 
Moreover, 
the {\pier additional} information $\Delta_{\Gamma } \mu _{\Gamma, \varepsilon} \in L^2(0,T;H_\Gamma )$ 
implies (see, e.g., \cite[p.~104]{Gri09})
$\mu _{\Gamma, \varepsilon} \in L^2 ( 0,T;H^2(\Gamma))$. 
Finally, this yields in 
particular that $\mu _{\Gamma, \varepsilon} \in L^2 ( 0,T;H^{3/2}(\Gamma))$, 
whence (quoting again \cite[Thm.~3.2, p.~1.79]{BG87}) we obtain 
$\mu _\varepsilon \in L^2 (0,T;H^2 (\Omega ))$, 
that is 
\begin{gather*}
	\mbox{\boldmath $ \mu $}_\varepsilon 
	\in L^2(0,T;\mbox{\boldmath $ W$}).
\end{gather*}
From the next subsection, we can {\pier proceed} with the a priori estimates. 

\subsection{A priori estimates}

In this subsection, we obtain the uniform estimates independent of {\pier $\varepsilon $}. 
We can {\pier adopt} the same strategy {\pier as in \cite{CC13, CF14, CF15, CGS14}. Here, we use systematically the relations 
\begin{gather}
\mbox{\boldmath $ u$}_\varepsilon = (u_\varepsilon, u_{\Gamma,\varepsilon  }) = 
\mbox{\boldmath $ v$}_\varepsilon    
	 +m_0 \mbox{\boldmath $ 1$} = (v_\varepsilon +m_0,v_{\Gamma,\varepsilon  }+m_0).
	 \nonumber
\end{gather}
}%

\paragraph{Lemma 4.1.}
{\it There exist{\pier s} a positive constant $M_1$, independent of $\varepsilon \in (0,1]$, such that}
\begin{align}
	& \varepsilon^{1/2}|\mbox{\boldmath $ v$}_\varepsilon 
	|_{L^\infty (0,T;\mbox{\scriptsize \boldmath $ H$}_0)}
	+ |\mbox{\boldmath $ v$}_\varepsilon|_{L^\infty (0,T;\mbox{\scriptsize \boldmath $ V$}_0^*)}
	+|\mbox{\boldmath $ v$}_\varepsilon |_{L^2(0,T;\mbox{\scriptsize \boldmath $ V$}_0)}
	\nonumber
	\\ 
	& {\pier{}+ \bigl| \beta _\varepsilon (u_\varepsilon ) 
	\bigr|_{{\colli L^1}(0,T;L^1(\Omega ))}
	+ \bigl| \beta _{\Gamma, \varepsilon} (u_{\Gamma, \varepsilon} ) 
	\bigr|_{{\colli L^1}(0,T;L^1(\Gamma))}}
	 \le M_1.
	 \label{pier6}
\end{align}

\paragraph{Proof.} We test \eqref{ape} at time $s \in (0,T)$ by 
$\mbox{\boldmath $ v$}_\varepsilon (s) \in \mbox{\boldmath $ V$}_0$. 
Then from {\colli \eqref{poin},} \eqref{proj} and the definition of the 
subdifferential with $\varphi ({\revis \mbox{\boldmath $ 0$}})=0$ we see that 
\begin{align}
	&{\colli {\revis \bigl (} \varepsilon \mbox{\boldmath $ v$}_\varepsilon ' (s ), 
	\mbox{\boldmath $ v$}_\varepsilon(s ) {\revis \bigr )}_{\mbox{\scriptsize \boldmath $ H$}_0} 
	+ {\revis \bigl (} \mbox{\boldmath $ v$}_\varepsilon'(s ), 
	\mbox{\boldmath $ v$}_\varepsilon(s ) {\revis \bigr )}_{\mbox{\scriptsize \boldmath $ V$}^*_0}}
	+ \varphi \bigl ( \mbox{\boldmath $ v$}_\varepsilon (s) \bigr )
	\nonumber \\
&	
	{\colli {}+ \bigl( \mbox{\boldmath $ \beta $}_\varepsilon 
	\bigl ( \mbox{\boldmath $ u$}_\varepsilon (s) \bigr ), 
	\mbox{\boldmath $ v$}_\varepsilon (s) \bigr)_{\! 
	\mbox{\scriptsize \boldmath $ H$}}}
	\le \bigl( \mbox{\boldmath $ f$}(s)-
	\mbox{\boldmath $ \pi $}
	\bigl( \mbox{\boldmath $ u$}_\varepsilon (s) \bigr), 
	\mbox{\boldmath $ v$}_\varepsilon (s) 
	\bigr)_{\! \mbox{\scriptsize \boldmath $ H$}} 
	\quad \mbox{for a.a.\  }s \in {\pier (0,T)}.
	\label{pier21}
\end{align} 
Now, let us recall the assumptions 
$D(\beta _\Gamma ) \subseteq D(\beta )$ of {\rm (A6)} {\pier and}
$m_0 \in {\rm int} D(\beta _\Gamma )$ of {\rm (A7)}. 
Hence, we can exploit the useful inequalities (a proof is given in  
\cite[Sect.~5]{GMS09}) stating the existence of two constants 
$\delta _0>0$ and $c_1>0$ such that 
\begin{gather*}
	\beta _\varepsilon (r)(r-m_0) \ge \delta _0 
	\bigl | \beta _\varepsilon(r) \bigr | -c_1,
	\quad 
	\beta _{\Gamma, \varepsilon} (r)(r-m_0) \ge \delta _0 
	\bigl | \beta _{\Gamma, \varepsilon} (r) \bigr | -c_1
\end{gather*}
for all $r \in \mathbb{R}$ and $\varepsilon \in (0,1]$. 
Therefore, we easily have
\begin{align}
\lefteqn{ 
	{\revis \bigl(} \mbox{\boldmath $ \beta $}_\varepsilon 
	\bigl ( \mbox{\boldmath $ u$}_\varepsilon (s) \bigr), 
	\mbox{\boldmath $ v$}_\varepsilon (s) {\revis \bigr )}_{\! \mbox{\scriptsize \boldmath $ H$}} 
} \nonumber \\
	& =  \int_{\Omega }^{} \beta _\varepsilon 
	\bigl ( u_\varepsilon (s) \bigr ) 
	(u_\varepsilon (s) - m_0) dx 
	+ \int_{\Gamma}^{} \beta _{\Gamma, \varepsilon} 
	\bigl ( u_{\Gamma, \varepsilon} (s) \bigr ) 
	 (u_{\Gamma, \varepsilon} (s) - m_0 ) d\Gamma 
	\nonumber \\
& \ge 
	\delta _0 \int_{\Omega }^{} \bigl| 
	 \beta _\varepsilon \bigl( u_\varepsilon(s) \bigr) 
	\bigr| dx -c_1 |\Omega |
	+ \delta _0 \int_{\Gamma }^{} \bigl| 
	 \beta _{\Gamma, \varepsilon} \bigl( u_{\Gamma, \varepsilon}(s) \bigr) 
	\bigr| d\Gamma  -c_1 |\Gamma |	\label{pier12} 
\end{align}
for a.a.\ $s \in {\pier (0,T)}$.
Moreover, {\pier in view of} the assumption {\rm (A4)} and {\pier the} compactness inequality \eqref{compin}, 
there exists a positive constant $\tilde{M}_1$, depending {\pier only} on $L$, $L_\Gamma $, 
$\pi(m_0)$, $\pi_\Gamma  (m_0)$, $|\Omega |$ and $|\Gamma|$, such that 
\begin{align}
\lefteqn{ 
	\bigl( \mbox{\boldmath $ f$}(s)-
	\mbox{\boldmath $ \pi $}
	\bigl( \mbox{\boldmath $ u$}_\varepsilon (s) \bigr), 
	\mbox{\boldmath $ v$}_\varepsilon (s) 
	\bigr)_{\! \mbox{\scriptsize \boldmath $ H$}} 
} \nonumber \\
	& \le \frac{1}{2}\int_{\Omega }^{} 
	\bigl | f (s) \bigr |^2 dx 
	+ \frac{1}{2} \int_{\Omega }^{} \bigl | v_\varepsilon (s) \bigr|^2 dx 
	+ L \int_{\Omega }^{}\bigl | v_\varepsilon (s) \bigr|^2 dx 
	+ \bigl|\pi (m_0)\bigr|  \int_{\Omega }^{}\bigl | v_\varepsilon (s) \bigr| dx 
	\nonumber \\
	& \quad {}
	+ \frac{1}{2}\int_{\Gamma }^{} 
	\bigl | f_\Gamma  (s) \bigr |^2 d\Gamma 
	+ \frac{1}{2} \int_{\Gamma}^{} \bigl | v_{\Gamma, \varepsilon }(s) \bigr|^2 d\Gamma  
	+ L_\Gamma  \int_{\Gamma }^{}\bigl | v_{\Gamma, \varepsilon } (s) \bigr|^2 d\Gamma 
	+ \bigl|\pi_\Gamma  (m_0) \bigr| \int_{\Gamma }^{}\bigl | v_{\Gamma, \varepsilon } (s) \bigr| d\Gamma  
	\nonumber \\
	& \le \tilde{M}_1
	\Bigl( 1 +  \bigl |\mbox{\boldmath $ f$}(s) 
	\bigr |_{\mbox{\scriptsize \boldmath $ H$}}^2 
	+ \bigl | \mbox{\boldmath $ v$}_{\varepsilon }(s) 
	\bigr |_{\mbox{\scriptsize \boldmath $ V$}_0^*}^2 \Bigr) 
	+ \frac{1}{4}\bigl | \mbox{\boldmath $ v$}_{\varepsilon }(s) 
	\bigr |_{\mbox{\scriptsize \boldmath $ V$}_0}^2 
	\label{pier22}
\end{align}
for a.a.\ $s \in {\pier (0,T)}$. 
Therefore,  collecting {\colli\eqref{pier21}--\eqref{pier22} it is straightforward to}
obtain 
\begin{align*}
	&\varepsilon \frac{d}{ds} \bigl 
	|\mbox{\boldmath $ v$}_\varepsilon (s ) \bigr |_{\mbox{\scriptsize \boldmath $ H$}_0}^2  
	+ \frac{d}{ds} \bigl |\mbox{\boldmath $ v$}_\varepsilon (s) 
	\bigr |_{\mbox{\scriptsize \boldmath $ V$}_0^*}^2 
	+ \frac{1}{2} \bigl | \mbox{\boldmath $ v$}_\varepsilon (s) 
	\bigr |_{\mbox{\scriptsize \boldmath $ V$}_0}^2
	\\
	&{\pier{}+ 2 \delta _0 \int_{\Omega }^{} {\revis \bigl |}
	 \beta _\varepsilon \bigl( u_\varepsilon(s) \bigr) 
	{\revis \bigr |} dx
	+ 2 \delta _0 \int_{\Gamma }^{} {\revis \bigl |}
	 \beta _{\Gamma, \varepsilon} \bigl( u_{\Gamma, \varepsilon}(s) \bigr) 
	{\revis \bigr |} d\Gamma  }
	\\
	&\le   {\pier {}2  c_1 {\revis \bigr (} |\Omega |  +  |\Gamma | {\revis \bigr )}} 
	+ 2 \tilde{M}_1
	\Bigl( 1 +  \bigl |\mbox{\boldmath $ f$}(s) 
	\bigr |_{\mbox{\scriptsize \boldmath $ H$}}^2 
	+ \bigl | \mbox{\boldmath $ v$}_{\varepsilon }(s) 
	\bigr |_{\mbox{\scriptsize \boldmath $ V$}_0^*}^2 \Bigr) 
	\quad \mbox{for a.a.\ }s \in {\pier (0,T)}.
\end{align*} 
Thus, the {G}ronwall inequality {\pier leads us} to the conclusion. \hfill $\Box$
\vspace{2mm}

Using the regularity assumption {\rm (A5)} for $\mbox{\boldmath $ f$}$, we {\pier can show} the following estimate. 

\paragraph{Lemma 4.2.}
{\it There exist{\pier s} a positive constant $M_2$, 
independent of $\varepsilon \in (0,1]$, such that}
\begin{align}
	&\varepsilon ^{1/2}|\mbox{\boldmath $ v$}_\varepsilon ' |_{L^2(0,T;\mbox{\scriptsize \boldmath $ H$}_0)} 
	+ |\mbox{\boldmath $ v$}_\varepsilon '| _{L^2(0,T;\mbox{\scriptsize \boldmath $ V$}_0^*)} 
	+ |\mbox{\boldmath $ v$}_\varepsilon |_{L^\infty (0,T;\mbox{\scriptsize \boldmath $ V$}_0)}  \nonumber
	\\ 
	&{}+ \bigl| \widehat{\beta} _\varepsilon (u_\varepsilon) 
	\bigr|_{L^\infty (0,T;L^1(\Omega ))} 
	+ \bigl| \widehat{\beta} _{\Gamma, \varepsilon} (u_{\Gamma, \varepsilon}) 
	\bigr|_{L^\infty (0,T;L^1(\Gamma ))} \le M_2.
\label{pier5}
\end{align}

\paragraph{Proof.} 
We test \eqref{ape} {\pier at} time $s \in (0,T)$ by 
$\mbox{\boldmath $ v$}_\varepsilon' (s) \in \mbox{\boldmath $ H$}_0$. 
Then, {\pier with the help of \eqref{proj}} we have 
\begin{align*}
	&\varepsilon  \bigl 
	|\mbox{\boldmath $ v$}_\varepsilon '(s ) 
	\bigr |_{\mbox{\scriptsize \boldmath $ H$}_0}^2  
	+  \bigl |\mbox{\boldmath $ v$}_\varepsilon '(s) 
	\bigr |_{\mbox{\scriptsize \boldmath $ V$}_0^*}^2 
	+ \frac{d}{ds} \varphi \bigl ( \mbox{\boldmath $ v$}_\varepsilon (s) \bigr )
	\\
	&\quad {}+\frac{d}{ds} \int_{\Omega }^{} 
	\widehat{\beta }_\varepsilon {\revis \bigl (}v_\varepsilon (s)+m_0{\revis \bigr )}dx
	+
	\frac{d}{ds} \int_{\Gamma }^{} 
	\widehat{\beta }_{\Gamma, \varepsilon} 
	{\revis \bigl (}v_{\Gamma, \varepsilon} (s)+m_0{\revis \bigr )}d\Gamma 
	\\
	&\le {\pier \frac{1}{4}} \bigl| \mbox{\boldmath $ v$}_\varepsilon '(s) 
	\bigr|_{\mbox{\scriptsize \boldmath $ V$}_0^*} ^2
	+ \int_{\Omega }^{} \bigl| \nabla \pi \bigl (v_\varepsilon (s)+m_0 \bigr) 
	\bigr|^2 dx 
	+ \int_{\Gamma }^{} \bigl| \nabla_\Gamma {\pier \pi_\Gamma} \bigl ({\pier v_{\Gamma, \varepsilon}} (s)+m_0\bigr) 
	\bigr|^2 d\Gamma  
	+ {\pier \bigl ( \mbox{\boldmath $ f$}(s),  \mbox{\boldmath $ v$}_\varepsilon '(s) 
	\bigr) _{\mbox{\scriptsize \boldmath $ H$}} }
\end{align*}
for a.a.\ $s \in {\pier (0,T)} $. 
Integrating it over $(0,t)$ with respect to $s$, {\pier in view  of 
\eqref{prim}, (A2) and (A7)} we infer that 
\begin{align}
	&\varepsilon \int_{0}^{t} \bigl 
	|\mbox{\boldmath $ v$}_\varepsilon '(s ) 
	\bigr |_{\mbox{\scriptsize \boldmath $ H$}_0}^2 ds  
	+ \frac{3}{4} \int_{0}^{t} \bigl |\mbox{\boldmath $ v$}_\varepsilon '(s) 
	\bigr |_{\mbox{\scriptsize \boldmath $ V$}_0^*}^2 ds 
	{\pier {}+ \frac12 {\revis \bigl |} 
	\mbox{\boldmath $ v$}_\varepsilon (t) 
	{\revis \bigr |}_{\mbox{\scriptsize \boldmath $ V$}_0}^2} 
	+ 
	\int_{\Omega }^{} 
	\widehat{\beta }_\varepsilon {\revis \bigl (}{\pier u_\varepsilon (t)}{\revis \bigr )}dx
	+
	\int_{\Gamma }^{} 
	\widehat{\beta }_{\Gamma, \varepsilon} {\revis \bigl (}
	{\pier u_{\Gamma, \varepsilon} (t)}{\revis \bigr )}d\Gamma 
	\nonumber 
	\\
	&\le 
	{\pier {} \frac12 \vert \mbox{\boldmath $ v$}_0 \vert_{\mbox{\scriptsize \boldmath $ V$}_0}^2} 
	+\int_{\Omega }^{} 
	\widehat{\beta }({\pier u_0} )dx
	+
	\int_{\Gamma }^{} 
	\widehat{\beta }_{\Gamma} ({\pier u_{0\Gamma}})d\Gamma 
	\nonumber
	\\
	&\quad {} + L^2 \int_{0}^{t} \! \! \int_{\Omega }^{} 
	\bigl| \nabla v_\varepsilon (s) \bigr|^2 dx ds
	+ L^2_\Gamma \int_{0}^{t} \! \! \int_{\Gamma }^{} 
	\bigl| \nabla_\Gamma v_\varepsilon (s) \bigr|^2 d\Gamma ds 
	+ \int_{0}^{t} {\pier \bigl ( \mbox{\boldmath $ f$}(s), 
	 \mbox{\boldmath $ v$}_\varepsilon '(s) 
	\bigr) _{\mbox{\scriptsize \boldmath $ H$}} } ds
	\label{pier15}
\end{align}
for all $t \in [0,T]$. 
{\pier Let us now recall (A5):  if  $\mbox{\boldmath $ f$} \in W^{1,1} (0,T; \mbox{\boldmath $ H$} )  $, 
we can  integrate the last term of \eqref{pier15} 
by parts in time.  With the help of 
the Young inequality and \eqref{poin} we easily obtain 
\begin{align*}
&\int_{0}^{t} {\pier \bigl ( \mbox{\boldmath $ f$}(s), 
	 \mbox{\boldmath $ v$}_\varepsilon '(s) 
	\bigr) _{\mbox{\scriptsize \boldmath $ H$}} } ds 
	\\
&= \bigl ( \mbox{\boldmath $ f$}(t), 
	 \mbox{\boldmath $ v$}_\varepsilon (t)  
	\bigr) _{\mbox{\scriptsize \boldmath $ H$}}  
- \bigl ( \mbox{\boldmath $ f$}(0), 
	 \mbox{\boldmath $ v$}_0 
	\bigr) _{\mbox{\scriptsize \boldmath $ H$}}  
- \int_{0}^{t} {\pier \bigl ( \mbox{\boldmath $ f$}'(s), 
	 \mbox{\boldmath $ v$}_\varepsilon (s) 
	\bigr) _{\mbox{\scriptsize \boldmath $ H $}} } ds
		\\
&\le 
\frac14  \bigl |	  \mbox{\boldmath $ v$}_\varepsilon (t)  
   \bigr |^2 _{\mbox{\scriptsize \boldmath $ V$}_0}  
+ \frac14  \bigl |	  \mbox{\boldmath $ v$}_0 \bigr |^2 _{\mbox{\scriptsize \boldmath $ H$}}  	
+ {\revis \left( \frac{1}{c_p} +1 \right)} \bigl | \mbox{\boldmath $ f$} \bigr |^2_{C ([0,T];  \mbox{\scriptsize \boldmath $ H$})	}
+  
\frac1{c_p^{\, 1/2}}\int_{0}^{t} {\pier \bigl | \mbox{\boldmath $ f$}'(s) \bigr |_{\mbox{\scriptsize \boldmath $ H$}}
	 \bigl | \mbox{\boldmath $ v$}_\varepsilon (s) 
	\bigr | _{\mbox{\scriptsize \boldmath $ V$}_0} } ds.
\end{align*}
Combinig this with \eqref{pier15} and applying the Gronwall lemma in the form of 
\cite[Lemme~A.5, p.~157]{Bre73}, there {\pier is} a positive constant $M_2$, depending {\pier only} on 
$|\mbox{\boldmath $ v$}_0|_{\mbox{\scriptsize \boldmath $ V$}_0}$, $c_p$,  
$|\widehat{\beta }(u_0)|_{L^1(\Omega )}$, 
$|\widehat{\beta }_\Gamma (u_{0\Gamma })|_{L^1(\Gamma )}$, 
$L$, $L_\Gamma $, $M_1$ and $|\mbox{\boldmath $ f$}|_{W^{1,1}(0,T;\mbox{\scriptsize \boldmath $ H$})}$,
such that {\pier \eqref{pier5} holds.} On the other hand, if 
 $ \mbox{\boldmath $ f$} \in L^2(0,T;\mbox{\boldmath $ V$}) $ in (A5)
 it suffices to observe that 
$$
\int_{0}^{t} {\pier \bigl ( \mbox{\boldmath $ f$}(s), 
	 \mbox{\boldmath $ v$}_\varepsilon '(s) 
	\bigr) _{\mbox{\scriptsize \boldmath $ H$}} } ds =
\int_{0}^{t} {\pier \bigl\langle 
	 \mbox{\boldmath $ v$}_\varepsilon '(s),
	 {\revis \mbox{\boldmath $ P$}\mbox{\boldmath $ f$}(s)} 
	\bigr \rangle _{\mbox{\scriptsize \boldmath $ V$}_0^*,\mbox{\scriptsize \boldmath $ V$}_0 }} ds 
	\le \frac{1}{4} \int_{0}^{t} \bigl |\mbox{\boldmath $ v$}_\varepsilon '(s) 
	\bigr |_{\mbox{\scriptsize \boldmath $ V$}_0^*}^2 ds + 
	 \bigl |\mbox{\boldmath $ f$}
	\bigr |_{L^2(0,T; \mbox{\scriptsize \boldmath $ V$})}^2 	
$$
and collecting this and \eqref{pier15} leads to the estimate \eqref{pier5} with a small change in the dependencies of $M_2$.}\hfill $\Box$ 

\paragraph{Lemma 4.3.}
{\it There exist {\pier two} positive constants ${\pier M_3}$ and ${\pier M_4}$,
independent of $\varepsilon \in (0,1]$, such that}
\begin{gather}
 {\colli{} \bigl| \beta _\varepsilon (u_\varepsilon ) 
	\bigr|_{{\colli L^2}(0,T;L^1(\Omega ))}
	+ \bigl| \beta _{\Gamma, \varepsilon} (u_{\Gamma, \varepsilon} ) 
	\bigr|_{{\colli L^2}(0,T;L^1(\Gamma))}
	 \le M_3,} \label{pier23}
   \\
	{\colli |\omega _\varepsilon |_{L^2(0,T)}  +  
	\bigl | \mbox{\boldmath $ \mu $}_\varepsilon 
	\bigr |_{L^2(0,T;\mbox{\scriptsize \boldmath $ V$} )} \le M_4}.
	\label{pier13}
\end{gather}
\paragraph{Proof.} {\colli Recalling \eqref{pier21} and \eqref{pier12}, we easily infer that 
\begin{align}
&\delta _0 \int_{\Omega }^{} \bigl| 
	 \beta _\varepsilon \bigl( u_\varepsilon(s) \bigr) 
	\bigr| dx
	+ \delta _0 \int_{\Gamma }^{} \bigl| 
	 \beta _{\Gamma, \varepsilon} \bigl( u_{\Gamma, \varepsilon}(s) \bigr) 
	\bigr| d\Gamma 
	\nonumber \\
&\le  c_1 {\revis \bigl(} |\Omega | +  |\Gamma |{\revis \bigr)} 	
	+ \bigl( \mbox{\boldmath $ f$}(s)- 	\mbox{\boldmath $ \pi $}
	\bigl( \mbox{\boldmath $ u$}_\varepsilon (s) \bigr)  - 
	\varepsilon \mbox{\boldmath $ v$}_\varepsilon ' (s )  , 
	\mbox{\boldmath $ v$}_\varepsilon (s) 
	\bigr)_{\! \mbox{\scriptsize \boldmath $ H$}} 
	-  ( \mbox{\boldmath $ v$}_\varepsilon'(s ), 	
	\mbox{\boldmath $ v$}_\varepsilon(s ) )_{\mbox{\scriptsize \boldmath $ V$}^*_0} 
	\label{pier24}
\end{align} 
for {\revis a.a.\ $s \in (0,T)$}. Hence, by squaring we have 
\begin{align*}
	& \left( \delta _0 \int_{\Omega }^{} {\revis \bigl| }
	 \beta _\varepsilon \bigl( u_\varepsilon(s) \bigr) 
	{\revis \bigr|} dx
	+ 
	\delta _0 \int_{\Gamma }^{} {\revis \bigl| }
	 \beta _{\Gamma, \varepsilon} \bigl( u_{\Gamma, \varepsilon}(s) \bigr) 
	{\revis \bigr|} d\Gamma \right)^2
	 \nonumber \\
	& \le 3 \/ c_1^{\,2} \bigl( |\Omega |+|\Gamma |\bigr)^2
	+ 
	 9 \/\Bigl(\bigl| \mbox{\boldmath $ f$}(s)
	\bigr| _{\mbox{\scriptsize \boldmath $ H$}}^2 
	+
	\bigl| \mbox{\boldmath $ \pi $}
	\bigl (\mbox{\boldmath $ u$}_\varepsilon (s) \bigr )
	\bigr|_{\mbox{\scriptsize \boldmath $ H$}}^2
	+\varepsilon^2
	\bigl| 
	\mbox{\boldmath $ v$}_\varepsilon' (s)
	\bigr|_{\mbox{\scriptsize \boldmath $ H$}_0}^2
	\Bigr)
	\bigl| 
	\mbox{\boldmath $ v$}_\varepsilon (s)
	\bigr|_{\mbox{\scriptsize \boldmath $ H$}_0}^2
	\\
	& \quad {}
	+
	3 \/\bigl| 
	\mbox{\boldmath $ v $}_\varepsilon' (s) 
	\bigr|_{\mbox{\scriptsize \boldmath $ V$}_0^*}^2
	\bigl| 
	\mbox{\boldmath $ v $}_\varepsilon (s)
	\bigr|_{\mbox{\scriptsize \boldmath $ V$}_0^*}^2
	\quad \mbox{for a.a.\ }s \in (0,T).  
\end{align*}
Then, due to Lemma~4.2, there exists a positive constant $M_3$, 
depending on $\delta _0$, $c_1$, $T$, $|\Omega |$, $|\Gamma |$, $M_2$, 
$L$, $L_\Gamma $, $|\pi (m_0)|$, $|\pi _\Gamma (m_0)|$, $|\mbox{\boldmath $ f$}|_{L^2(0,T;\mbox{\scriptsize \boldmath $ H$})}$ and 
independent of $\varepsilon \in (0,1]$, such that \eqref{pier23} holds. Next,
from the definition of $\omega _\varepsilon $, give{\pier n} by \eqref{ome}, 
we have 
\begin{multline*}
	{\pier |\omega_\varepsilon(t)  |^2
   \le 
	\frac{6}{ {\revis \bigl( } |\Omega |+|\Gamma |{\revis \bigr) }^2 }
	\left\{ 
	\bigl\vert  \beta _\varepsilon \bigl( u_\varepsilon (t) \bigr ) \bigr\vert^2_{L^1(\Omega )} 
	+ \bigl\vert  \pi \bigl( u_\varepsilon (t) \bigr ) 	\bigr\vert^2_{L^1(\Omega )}
	+
	|\Omega |  \bigl|  	f (t) 	\bigr|^2_{H} 
    \right. }
	\\
	 {\pier \left.  	{}+
	 \bigl|  
	\beta _{\Gamma, \varepsilon} \bigl( u_{\Gamma, \varepsilon} (t) \bigr )
	\bigr|^2_{L^1(\Gamma)} +  \bigl|  
	\pi_\Gamma  \bigl( u_{\Gamma ,\varepsilon} (t) \bigr )
	\bigr|^2_{L^1(\Gamma  )} 
	+
	|\Gamma | \bigl|  
	f_\Gamma  (t) 
	\bigr|^2_{H_\Gamma}
	\right\}.  }
\end{multline*} 
Then, by integrating over $(0,T)$, it follows that there is a positive constant ${\tilde M_4}$, 
depending only on $|\Omega |$, $|\Gamma |$, $M_2$, $M_3$,   
$L$, $L_\Gamma $, $|\pi (m_0)|$, $|\pi _\Gamma (m_0)|$ and $|\mbox{\boldmath $ f$}|_{L^2(0,T;\mbox{\scriptsize \boldmath $ H$})}$, such that 
$$ |\omega _\varepsilon |_{L^2(0,T)} \le {\tilde M_4}.$$
At this point, we test \eqref{eonvs} at time $s \in (0,T)$ by 
$\mbox{\boldmath $ P$} 
\mbox{\boldmath $ \mu $}_\varepsilon (s) \in \mbox{\boldmath $ V$}_0$
and obtain 
\begin{gather*}
	\bigl \langle 
	\mbox{\boldmath $ F$} \bigl( 
	\mbox{\boldmath $ P$} \mbox{\boldmath $ \mu $}_\varepsilon (s) 
	\bigr),
	\mbox{\boldmath $ P$} \mbox{\boldmath $ \mu $}_\varepsilon (s)
	\bigr \rangle_{\mbox{\scriptsize \boldmath $ V$}_0^*,
	\mbox{\scriptsize \boldmath $ V$}_0}
	= - \bigl \langle 
	\mbox{\boldmath $ v $}_\varepsilon' (s) ,
	\mbox{\boldmath $ P$} \mbox{\boldmath $ \mu $}_\varepsilon (s)
	\bigr \rangle_{\mbox{\scriptsize \boldmath $ V$}_0^*,
	\mbox{\scriptsize \boldmath $ V$}_0}
\end{gather*}
for {\revis a.a.\ $s\in (0,T)$}. Then, in view of \eqref{poin} we deduce that 
\begin{gather*}
	\bigl | 
	\mbox{\boldmath $ P$} \mbox{\boldmath $ \mu $}_\varepsilon (s)
	\bigr |_{\mbox{\scriptsize \boldmath $ V$}_0}
	\le  \bigl |
	\mbox{\boldmath $ v $}_\varepsilon' (s) 
	\bigr |_{\mbox{\scriptsize \boldmath $ V$}_0^*}
\end{gather*}
and 
\begin{align*}
	\bigl | \mbox{\boldmath $ \mu $}_\varepsilon (s)
	\bigr |_{\mbox{\scriptsize \boldmath $ V$} } 
	& \le \bigl | \mbox{\boldmath $ P$} \mbox{\boldmath $ \mu $}_\varepsilon (s)
	\bigr |_{\mbox{\scriptsize \boldmath $ V$} } 
	+ \bigl| \omega_\varepsilon (s) \mbox{\boldmath $ 1$}
	\bigr |_{\mbox{\scriptsize \boldmath $ V$} } 
	\\ 
	& \le 
	c_p^{-1/2} 	\bigl | \mbox{\boldmath $ P$} \mbox{\boldmath $ \mu $}_\varepsilon (s)
	\bigr |_{\mbox{\scriptsize \boldmath $ V$}_0 } 
	+ \bigl( |\Omega | +|\Gamma | \bigr)^{1/2} |\omega _\varepsilon (s) | 
\end{align*} 
{\pier whence \eqref{pier13}  follows by squaring and integrating over $(0,T) $,
on account of the estimate for $|\mbox{\boldmath $ v$}_\varepsilon '| _{L^2(0,T;\mbox{\scriptsize \boldmath $ V$}_0^*)}$ in \eqref{pier5}.}}
\hfill $\Box$

\paragraph{Lemma 4.4.}
{\it There exist {\pier two} positive constants ${\pier M_5}$ and ${\pier M_6}$,
independent of $\varepsilon \in (0,1]$, such that}
\begin{gather}
	\bigl| \beta _\varepsilon (u_\varepsilon) \bigr |_{L^2(0,T;H)} 
	+ \bigl |\beta _\varepsilon (u_{\Gamma,\varepsilon }) \bigr |_{L^2(0,T;H_\Gamma )} 
	\le {\pier M_5}, \label{pier7}
	\\
	{\pier |\Delta v_\varepsilon |_{L^2(0,T;H)}} + 
	|v_\varepsilon |_{L^2(0,T;H^{3/2}(\Omega ))}
	+ 
	|\partial_\nu v_\varepsilon |_{L^2(0,T;H_\Gamma )} 
	\le {\pier M_6}. \label{pier8}
\end{gather}

\paragraph{Proof.} We test \eqref{d5e} by 
$\beta _\varepsilon (u_\varepsilon) \in L^2(0,T;V)$ 
and {\pier use} \eqref{d7e},  {\pier by noting that
$(\beta _\varepsilon (u_\varepsilon ))_{|_\Gamma }
=\beta _\varepsilon (u_{\Gamma, \varepsilon}) \in L^2(0,T;V_\Gamma)$ as well.
Then, by integrating over $\Omega$, 
we infer that
\begin{align*}
	\lefteqn{ 
	\int_{\Omega }^{} \beta '_\varepsilon 
	\bigl( u_\varepsilon (s) \bigr)
	\bigl | \nabla v_\varepsilon (s) \bigr |^2 dx
	+ 
	\bigl | \beta _\varepsilon 
	\bigl (u_\varepsilon (s)\bigr) \bigr |^2_{H} 
	} \nonumber \\
	&
	+ \int_{\Gamma }^{} \beta '_\varepsilon 
	\bigl( u_{\Gamma, \varepsilon} (s) \bigr)
	\bigl | \nabla_\Gamma  v_{\Gamma,\varepsilon} (s) \bigr |^2 d\Gamma 
	+ 
	\int_{\Gamma }^{} 
	\beta _{\Gamma, \varepsilon } \bigl( u_{\Gamma, \varepsilon }(s)  \bigr)
	\beta _\varepsilon \big( u_{\Gamma, \varepsilon } (s) \bigr)
	d\Gamma 
	\nonumber \\
	& \le 
	 \bigl ( f(s ) + \mu _\varepsilon (s)  
	 - \varepsilon\partial_t v_\varepsilon (s)
	 - \pi \bigl (u_\varepsilon (s )
	\bigr ) ,\beta _\varepsilon 
	\bigl( u_\varepsilon (s ) \bigr) 
	\bigr )_{\! H} 
	\\
	& \quad  {}
	+
	\bigl ( f_\Gamma (s )+ \mu _{\Gamma ,\varepsilon }(s)
	- \varepsilon \partial_t v_{\Gamma ,\varepsilon }(s)	- \pi _\Gamma 
	\bigl (u_{\Gamma , \varepsilon} (s ) 
	\bigr ) ,
	\beta _\varepsilon 
	\bigl( u_{\Gamma, \varepsilon}(s )
	\bigr) 
	\bigr )_{\! H_\Gamma } 
\end{align*} 
for {\revis a.a.\ $s\in (0,T)$}. 
{\pier Recalling now the condition} \eqref{a4e} we deduce that 
\begin{align*} 
	&\int_{\Gamma }^{} 
	\beta _{\Gamma, \varepsilon } \bigl( u_{\Gamma, \varepsilon }(s) \bigr)
	\beta _\varepsilon \big( u_{\Gamma, \varepsilon } (s)\bigr)
	d\Gamma 
	 = 
	\int_{\Gamma }^{} 
	\bigl|
	\beta _{\Gamma, \varepsilon } \bigl( u_{\Gamma, \varepsilon }(s) \bigr)
	\bigr| 
	\bigl| 
	\beta _\varepsilon \big( u_{\Gamma, \varepsilon } (s) \bigr) 
	\bigr| 
	d\Gamma 
	\\ &\qquad{}
	 \ge
	\frac{1}{\varrho } \int_{\Gamma }^{} 
	\bigl | \beta _\varepsilon \bigr( u_{\Gamma ,\varepsilon } (s)  \bigr) \bigr | ^2 
	d\Gamma
	- \frac{c_0}{\varrho }\int_{\Gamma }^{} 
	\bigl |\beta _\varepsilon \bigr( u_{\Gamma ,\varepsilon } (s) \bigr) \bigr| 
	d\Gamma
	\\ &\qquad{}
	 \ge
	\frac{1}{2\varrho } \int_{\Gamma }^{} 
	\bigl | \beta _\varepsilon \bigr( u_{\Gamma ,\varepsilon } (s)  \bigr) \bigr | ^2 
	d\Gamma 
	- \frac{c_0^2}{2\varrho } \, |\Gamma | {\pier ,}
\end{align*}
because $\beta _\varepsilon (r)$ and $\beta _{\Gamma, \varepsilon }(r)$ have the same sign 
for all $r \in \mathbb{R}$. {Also, we observe} that 
\begin{gather*}
	\int_{\Omega }^{} \beta '_\varepsilon 
	\bigl( u_\varepsilon (s) \bigr)
	\bigl | \nabla v_\varepsilon (s) \bigr |^2 dx \ge 0,
	\quad \
	 \int_{\Gamma }^{} \beta '_\varepsilon 
	\bigl( u_{\Gamma, \varepsilon} (s) \bigr)
	\bigl | \nabla_\Gamma v_{\Gamma,\varepsilon} (s) \bigr |^2 d\Gamma  \ge 0 .
\end{gather*}
Moreover, by the Young inequality and the Lipschitz continuity of $\pi$ and $\pi_\Gamma$ we can find a positive constant $\tilde{M}_5$, 
independent of $\varepsilon \in (0,1]$, such that
\begin{align*} 
	\lefteqn{ 
	\bigl ( f(s ) + \mu_\varepsilon (s) 
	- \pi \bigl (u_\varepsilon (s )
	\bigr ) ,\beta _\varepsilon 
	\bigl( u_\varepsilon (s ) \bigr) 
	\bigr )_{\! H} 
	} \nonumber \\
	& \le  
	\frac{1}{2} 
	\bigl| \beta _\varepsilon 
	\bigl( u_\varepsilon (s )\bigr) 
	\bigr|_{H}^2 
	+
	\tilde{M}_5
	{\revis \bigl(} 
	1+{\revis \bigl|}f(s ){\revis \bigr|}_H^2 
	+ {\revis \bigl|}\mu _\varepsilon (s ){\revis \bigr|}_H^2 
	+ {\revis \bigl|}v_\varepsilon(s ) {\revis \bigr|}_H^2 
	+ \varepsilon {\revis \bigl|}\partial_t v_\varepsilon (s ) {\revis \bigr|}_H^2 
	{\revis \bigr)}, 
	\\[0.3cm]
	&\bigl ( f_\Gamma (s )+\mu _{\Gamma, \varepsilon} 
    {\pier (s)} 
	-\pi _\Gamma 
	\bigl (u_{\Gamma , \varepsilon} (s )
	\bigr ),
	\beta _\varepsilon 
	\bigl( u_{\Gamma, \varepsilon}(s)
	\bigr) 
	\bigr )_{\! H_\Gamma }
	\nonumber 
	\\
	& \le 
	\frac{1 }{4\varrho}
	\bigl| \beta _\varepsilon 
	\bigl( u_{\Gamma, \varepsilon }(s ) \bigr) 
	\bigr|_{H_\Gamma }^2  + 
	 \varrho \tilde{M}_5
	{\revis \bigl(} 
	1+{\revis \bigl|}f_\Gamma (s) {\revis \bigr|}_{H_\Gamma}^2 
	+ {\revis \bigl|}\mu _{\Gamma,\varepsilon }(s)  {\revis \bigr|}_{H_\Gamma}^2 
	+ {\revis \bigl|}v_{\Gamma,\varepsilon}(s) {\revis \bigr|}_{H_\Gamma}^2
	+ \varepsilon {\revis \bigl|} \partial_t v_{\Gamma,\varepsilon}(s) {\revis \bigr|}_{H_\Gamma}^2
	{\revis \bigr)} .
\end{align*} 
Combining the above inequalities, integrating in $(0,T)$ with respect to $s$
and recalling \eqref{pier5} and \eqref{pier13} we infer that
\begin{align*}
    \frac{1}{2}\int_{0}^{T} 
	\bigl | \beta _\varepsilon 
	\bigl (u_\varepsilon (s)\bigr) \bigr |^2_{H} ds 
	+ 
	 \frac{1}{4\varrho }
	\int_{0}^{T} 
	\bigl | \beta _\varepsilon \bigl( u_{\Gamma ,\varepsilon } (s)
	 \bigr) \bigr |_{H_\Gamma }^2 
	 ds
\end{align*} 
is bounded independently of  $\varepsilon$, that is},
there is a positive constant ${\pier M_5}$, 
independent of $\varepsilon \in (0,1]$, such that 
\begin{gather*}
	\bigl |\beta _\varepsilon (u_\varepsilon )\bigr |_{L^2(0,T;H)} 
	+ \bigl |\beta _\varepsilon (u_{\Gamma,\varepsilon })\bigr |_{L^2(0,T;H_\Gamma )} 
	\le {\pier M_5}.
\end{gather*}
Now, we can compare the terms in \eqref{d5e} and 
conclude that {\pier $
	|\Delta v_\varepsilon |_{L^2(0,T;H)}
$ is uniformly bounded as well.
Hence, with the help of Lemma~4.1 and applying the theory of the elliptic regularity 
(see, e.g., \cite[{\pier Thm.~3.2, p.~1.79}]{BG87}), we have that 
\begin{gather*}
	|v_\varepsilon |_{L^2(0,T;H^{3/2}(\Omega ))}
	\le {\pier \tilde M_6}
\end{gather*}
and, owing to the trace theory (see, e.g., \cite[{\pier Thm.~2.25, p.~1.62}]{BG87}), that
\begin{gather*}
	|\partial_\nu v_\varepsilon |_{L^2(0,T;H_\Gamma )} 
	\le {\pier \tilde M_6}
\end{gather*}
for some positive constant ${\pier \tilde M_6}$ that is independent of $\varepsilon \in (0,1]$. Thus, the lemma is completely proved.}
\hfill $\Box$

\paragraph{Lemma 4.5.}
{\it There exist {\pier two} positive constants ${\pier M_7}$ and ${\pier M_8}$,
independent of $\varepsilon \in (0,1]$, such that} 
\begin{gather}
	\bigl |\beta _{\Gamma, \varepsilon} (u_{\Gamma,\varepsilon }) 
	\bigr |_{L^2(0,T;H_\Gamma )} \le {\pier M_7},
    \quad \ 	
    |\mbox{\boldmath $ v$}_{\varepsilon} |_{L^2(0,T;\mbox{\scriptsize 
    \boldmath $ W$})}
	\le {\pier M_8}. \label{pier18}
\end{gather}

\paragraph{Proof.} We test \eqref{d7e} by 
$\beta _{\Gamma, \varepsilon} (u_{\Gamma, \varepsilon}) \in L^2(0,T;V_\Gamma)$ and 
integrate {\pier it} on the boundary {\pier $\Gamma$, deducing that
\begin{align*}
	& 
	\int_{\Gamma}^{} \beta '_{\Gamma, \varepsilon }
	\bigl( u_{\Gamma, \varepsilon} (s) \bigr)
	\bigl | \nabla_\Gamma v_{\Gamma, \varepsilon} (s) \bigr |^2 d\Gamma
	{\pier {}+ \bigl| \beta _{\Gamma, \varepsilon }
	\bigl( u_{\Gamma, \varepsilon}(s)
	\bigr) \bigr|_{H_\Gamma }^2  }\nonumber \\
	& \le
	\bigl ( f_\Gamma (s )+ \mu _\Gamma (s)-\varepsilon \partial_t v_{\Gamma , \varepsilon}(s)
	-  \partial_\nu  v_\varepsilon(s)-\pi _\Gamma 
	\bigl (u_{\Gamma , \varepsilon} (s ) 
	\bigr ),
	\beta _{\Gamma, \varepsilon }
	\bigl( u_{\Gamma, \varepsilon}(s ) 
	\bigr) 
	\bigr )_{\! H_\Gamma }\\
	& \le 
	\frac{1}{2}\bigl| \beta _{\Gamma, \varepsilon }
	\bigl( u_{\Gamma, \varepsilon}(s)
	\bigr) \bigr|_{H_\Gamma }^2 
	\\ &\quad {}
	+ 3 \left(
	\bigl |f_\Gamma (s ) \bigr|^2_{H_\Gamma } 
	+  \bigl| \mu _\Gamma (s) \bigr|_{H_\Gamma } ^2
    + \varepsilon 
	\bigl| \partial_t v_{\Gamma , \varepsilon}(s) \bigr|_{H_\Gamma }^2 
	+ 
	\bigl| \partial_\nu v_\varepsilon(s) \bigr|_{H_\Gamma }^2 
	+	\bigl| \pi _\Gamma 
	\bigl(u_{\Gamma , \varepsilon} (s ) \bigr) 
	\bigr|_{H_\Gamma}^2 \right)
\end{align*} 
for {\revis a.a.\ $s\in (0,T)$}. {\pier By {\revis neglecting} the first positive contribution and integrating over $(0,T)$, we find out that there {\pier is} a positive constant ${\pier M_7}$, 
depending only  on 
$|f_\Gamma |_{L^2(0,T;H_\Gamma )}$, ${\pier M_4}$, ${\pier M_2}$,} ${\pier M_6}$, $L_\Gamma $, 
$|\pi _\Gamma (m_0)|$, $|\Gamma |$ and $T$, such that} 
\begin{gather*}
	\bigl |\beta _{\Gamma, \varepsilon} (u_{\Gamma,\varepsilon }) 
	\bigr |_{L^2(0,T;H_\Gamma )} \le {\pier M_7}.
\end{gather*}
{\pier Thanks to \eqref{pier6} and \eqref{pier8}, by comparison in \eqref{d7e} we also infer that
$
	|\Delta_\Gamma v_{\Gamma, \varepsilon} |_{L^2(0,T;H_\Gamma )} 
$
is bounded independently of $\varepsilon \in (0,1]$, and 
	 consequently (see, e.g., \cite[Sect.~4.2]{Gri09}),
\begin{gather} 
	|v_{\Gamma, \varepsilon }|_{L^2(0,T;H^2(\Gamma ))} 
	 \le  \left( |v_{\Gamma, \varepsilon }|_{L^2(0,T;V_\Gamma)}^2+
	|\Delta _\Gamma v_{\Gamma, \varepsilon }|_{L^2(0,T;H_\Gamma)}^2 \right)^{1/2} 
	 \le \tilde{M}_8. \label{pier9}
\end{gather} 
for some constant $\tilde{M}_8.$ Then, using the theory of the 
elliptic regularity (see, e.g., \cite[{\pier Thm.~3.2, p.~1.79}]{BG87}, 
by virtue of \eqref{pier6} and \eqref{pier8} it follows that
$ {\pier  |v_\varepsilon |_{L^2(0,T;H^2(\Omega ))}} $
is uniformly bounded, whence (cf.~\eqref{pier9})
\begin{gather*}
	|\mbox{\boldmath $ v$}_\varepsilon 
	|_{L^2(0,T;\mbox{\scriptsize \boldmath $ W$})} \le {\pier M_8}
\end{gather*}
for some positive constant ${\pier M_8}$ independent of $\varepsilon \in (0,1]$. 
Hence, \eqref{pier18} is proved.} 
\hfill $\Box$

\subsection{Passage to the limit as $\varepsilon \to 0$}

In this subsection, we conclude the existence proof by {\pier passing to the limit in the approximating problem} as $\varepsilon \to 0$. Indeed, 
owing to the estimates stated {\pier in} Lemmas from 4.1 to 4.5, there 
exist a subsequence of $\varepsilon $ (not relabeled) and some 
limit functions {\pier $ \mbox{\boldmath $ v$}$, $\mbox{\boldmath $ \mu $}$,} 
$\xi$, $\xi _\Gamma$ and $\omega$ such that 
\begin{gather} 
	\mbox{\boldmath $ v$}_\varepsilon \to \mbox{\boldmath $ v$} \quad \mbox{weakly star in } 
	H^1(0,T;\mbox{\boldmath $ V$}_0^*) 
	\cap 
	L^\infty (0,T;\mbox{\boldmath $ V$}_0) 
	\cap 
	L^2 \bigl (0,T;\mbox{\boldmath $ W$} \bigr), 
	\label{wcv}
	\\
	\varepsilon \mbox{\boldmath $ v$}_{\varepsilon} \to  {\pier \mbox{\boldmath $ 0$}}
	\quad \mbox{strongly in } 
	H^1(0,T;\mbox{\boldmath $ H$}_0),
	\label{scev}
	\\
	\mbox{\boldmath $ \mu $}_\varepsilon \to \mbox{\boldmath $ \mu $} 
	\quad \mbox{weakly in } 
	L^2(0,T;\mbox{\boldmath $ V$}), 
	\label{wcmu}
	\\
	{\pier \beta_\varepsilon(u_\varepsilon)} \to \xi \quad \mbox{weakly in } 
	L^2(0,T;H),
	\label{wcb}
	\\ 
	{\pier \beta _{\Gamma, \varepsilon} (u_{\Gamma,\varepsilon }) }\to \xi _\Gamma 
	\quad \mbox{weakly in } 
	L^2(0,T;H_\Gamma ),
	\label{wcbg}
	\\
	\omega _\varepsilon \to \omega  \quad \mbox{weakly in } 
	L^2(0,T).
	\label{wco}
\end{gather} 
{\pier From \eqref{wcv}, due a well-known} compactness results 
(see, e.g., \cite[{\pier Sect.~8, Cor.~4}]{Sim87}) 
we obtain
\begin{gather} 
	\mbox{\boldmath $ v$}_\varepsilon \to \mbox{\boldmath $ v$} \quad \mbox{strongly in } 
	C\bigl( [0,T];\mbox{\boldmath $ H$}_0 \bigr) \cap L^2 (0,T;\mbox{\boldmath $ V$}_0),
	\label{scv}
\end{gather}
{\pier which also entails} 
\begin{gather} 
	\mbox{\boldmath $ u$}_\varepsilon \to \mbox{\boldmath $ u$} 
	=\mbox{\boldmath $ v$}+m_0 \mbox{\boldmath $ 1$}
	\quad \mbox{strongly in } 
	C\bigl( [0,T];\mbox{\boldmath $ H$} \bigr) \cap L^2 (0,T;\mbox{\boldmath $ V$}),
	\label{scu}
\end{gather}
as $\varepsilon \to 0$. 
We point out that \eqref{scv} implies {\pier $\mbox{\boldmath $ v$} (0) = \mbox{\boldmath $ v$}_0$, that is,}
\begin{gather*}
	v(0)=v_0 \quad \mbox{a.e.\ in } \Omega, 
	\quad 
	v_\Gamma (0) = v_{0\Gamma } 
	\quad \mbox{a.e.\ on } \Gamma.
\end{gather*}
Moreover, \eqref{scu} and 
the {L}ipschitz continuity of $\pi $ and $\pi _\Gamma $  {\pier ensure that}
\begin{gather*}
	\mbox{\boldmath $ \pi$} (\mbox{\boldmath $ u$}_\varepsilon) 
	\to \mbox{\boldmath $ \pi $}(\mbox{\boldmath $ u$}) \quad \mbox{strongly in }
	C \bigl( [0,T];\mbox{\boldmath $ H$} \bigr),
\end{gather*} 
{\pier whence (cf.~\eqref{ome}) $\omega =m(\mbox{\boldmath $ \xi $}
+\mbox{\boldmath $ \pi $}(\mbox{\boldmath $ u$})-\mbox{\boldmath $ f$})$ with $
\mbox{\boldmath $ \xi $} = (\xi, \xi_\gamma).$}
Moreover, by applying \cite[{\pier Prop.~2.2, p.~38}]{Bar10} and using 
\eqref{wcb}--\eqref{wcbg} with \eqref{scu}, we {\pier deduce that
\begin{gather*}
	\xi \in \beta (u)
	\quad \mbox{a.e.\ in } Q, \quad 
	\xi_\Gamma \in \beta_{\Gamma } (u_{\Gamma })
	\quad \mbox{a.e.\ on } \Sigma
\end{gather*}
due to the maximal monotonicity of $\beta$ and $\beta_\Gamma.$}
At this point, we can pass 
to the limit in {\pier \eqref{wf1e}--\eqref{wf2e}} obtaining 
\eqref{d1}--\eqref{d2}. 
Thus, it turns out that {\pier the triplet $(\mbox{\boldmath $ v$}, \mbox{\boldmath $ \mu  $}, \mbox{\boldmath $ \xi$})$} is a weak solution of {\rm {\pier (P)}}.

\subsection{Regularity result}

In this subsection, we try to prove a regularity 
estimate allowing us to {\pier fix a strong solution of {\rm {\pier (P)}}. 
Then, our third theorem deals with the existence of the strong solution. 
To this aim, we introduce} the following additional regularity assumptions for 
$\mbox{\boldmath $ f$}$ and $\mbox{\boldmath $ u$}_0$: 
\begin{itemize}
 \item[(A8)] $\mbox{\boldmath $ f$} \in H^1(0,T;\mbox{\boldmath $ H$})$; 
 \item[(A9)] 
$\mbox{\boldmath $ u$}_0 \in \mbox{\boldmath $ W$}$ and 
the family $\{ -\partial \varphi (\mbox{\boldmath $ u$}_0-m_0\mbox{\boldmath $ 1$})
-\mbox{\boldmath $ \beta $}_\varepsilon (\mbox{\boldmath $ u$}_0)-
\mbox{\boldmath $ \pi $}(\mbox{\boldmath $ u$}_0)+\mbox{\boldmath $ f$}(0) \, : \ {\pier \varepsilon \in (0,\varepsilon_0 ]\/}\}$ 
is bounded in ${\pier \mbox{\boldmath $ V$}}$ for some $\varepsilon _0 \in (0,1]$. 
\end{itemize} 
These assumptions, {\pier in particular (A9), can be compared with the 
analogous ones in \cite{CGS14}. Let us point out that (A8) entails the validity of (A1) and (A5).} 

\paragraph{Theorem 2.3.} 
{\it Under the assumptions {\rm (A3)}, {\rm (A4)}, {\rm (A6)}--{\rm (A9)}, 
there exist{\pier s} a strong solution of {\rm {\pier (P)}}. 
}

\paragraph{Proof.} Recall \eqref{ape}--\eqref{ice}
at times $s, s+h\in (0,T)$. Take {\pier the} difference of them and 
test the resultant by 
$\mbox{\boldmath $ v$}_\varepsilon(s+h) -\mbox{\boldmath $ v$}_\varepsilon (s)$ where 
$0<h<T-s$. 
Then, by virtue of {\rm (A4)} and the compactness inequality \eqref{compin}, {\pier we have that}
\begin{align*}
\lefteqn{ 
	\frac{\varepsilon }{2} \frac{d}{ds}\bigl 
	|\mbox{\boldmath $ v$}_\varepsilon (s+h)-\mbox{\boldmath $ v$}_\varepsilon (s) \bigr |_{\mbox{\scriptsize \boldmath $ H$}_0}^2  
	+ \frac{1}{2} \frac{d}{ds} \bigl |\mbox{\boldmath $ v$}_\varepsilon (s+h) 
	-\mbox{\boldmath $ v$}_\varepsilon (s) 
	\bigr |_{\mbox{\scriptsize \boldmath $ V$}_0^*}^2 
} \nonumber \\
	& \quad {}+ \bigl | 
	\mbox{\boldmath $ v$}_\varepsilon (s+h)
	-\mbox{\boldmath $ v$}_\varepsilon (s) \bigr |_{\mbox{\scriptsize \boldmath $ V$}_0}^2
	+ \bigl( \mbox{\boldmath $ \beta $}_\varepsilon 
	\bigl ( \mbox{\boldmath $ u$}_\varepsilon (s) \bigr )
	- \mbox{\boldmath $ \beta $}_\varepsilon 
	\bigl ( \mbox{\boldmath $ u$}_\varepsilon (s+h) \bigr ), 
	\mbox{\boldmath $ v$}_\varepsilon (s+h) - 
	\mbox{\boldmath $ v$}_\varepsilon (s)\bigr)_{\! \mbox{\scriptsize \boldmath $ H$}} 
	\\
	& \le (L+L_\Gamma ) \bigl | 
	\mbox{\boldmath $ v$}_\varepsilon (s+h)
	-\mbox{\boldmath $ v$}_\varepsilon (s) \bigr |_{\mbox{\scriptsize \boldmath $ H$}}^2
	+ \frac{1}{2}\bigl | 
	\mbox{\boldmath $ f$}_\varepsilon (s+h)
	-\mbox{\boldmath $ f$}_\varepsilon (s) \bigr |_{\mbox{\scriptsize \boldmath $ H$}}^2 
	+\frac{1}{2} \bigl | 
	\mbox{\boldmath $ v$}_\varepsilon (s+h)
	-\mbox{\boldmath $ v$}_\varepsilon (s) \bigr |_{\mbox{\scriptsize \boldmath $ H$}}^2 
	\\
	& \le 
	\frac{1}{2} \bigl | 
	\mbox{\boldmath $ f$}_\varepsilon (s+h)
	-\mbox{\boldmath $ f$}_\varepsilon (s) \bigr |_{\mbox{\scriptsize \boldmath $ H$}}^2 
	+
	\tilde{\delta } 
	\bigl | 
	\mbox{\boldmath $ v$}_\varepsilon (s+h)
	-\mbox{\boldmath $ v$}_\varepsilon (s) \bigr |_{\mbox{\scriptsize \boldmath $ V$}_0}^2 
	+c_{\tilde{\delta }} \bigl | 
	\mbox{\boldmath $ v$}_\varepsilon (s+h)
	-\mbox{\boldmath $ v$}_\varepsilon (s) \bigr |_{\mbox{\scriptsize \boldmath $ V$}_0^*}^2
\end{align*} 
for a.a.\ $s \in {\pier (0,T)}$, where $\tilde{\delta }>0$ and 
$c_{\tilde{\delta }}>0$ is a constant depending on $L$, $L_\Gamma $. 
Taking $\tilde{\delta }=1/2$, dividing the resultant by $h^2$, 
integrating it over $(0,t)$ with respect to $s$, 
and using the monotonicity of $\mbox{\boldmath $ \beta $}$, we infer that 
\begin{align*}
	&
	\frac{\varepsilon}{2}  \left| 
	\frac{\mbox{\boldmath $ v$}_\varepsilon (t+h)-\mbox{\boldmath $ v$}_\varepsilon (t)}{h} 
	\right|_{\mbox{\scriptsize \boldmath $ H$}_0}^2  
	+ \frac{1}{2}\left| \frac{\mbox{\boldmath $ v$}_\varepsilon (t+h) 
	-\mbox{\boldmath $ v$}_\varepsilon (t) }{h}
	\right |_{\mbox{\scriptsize \boldmath $ V$}_0^*}^2 
	+ \frac{1}{2}\int_{0}^{t}\left| \frac{\mbox{\boldmath $ v$}_\varepsilon (s+h) 
	-\mbox{\boldmath $ v$}_\varepsilon (s) }{h}
	\right |_{\mbox{\scriptsize \boldmath $ V$}_0}^2 ds 
	\\
	& \le 
	\frac{\varepsilon}{2}  \left| 
	\frac{\mbox{\boldmath $ v$}_\varepsilon (h)-\mbox{\boldmath $ v$}_{0}}{h} 
	\right|_{\mbox{\scriptsize \boldmath $ H$}_0}^2 
	+ \frac{1}{2}  \left| 
	\frac{\mbox{\boldmath $ v$}_\varepsilon (h)-\mbox{\boldmath $ v$}_{0}}{h} 
	\right|_{\mbox{\scriptsize \boldmath $ V$}_0^*}^2 
		\\ 
	& \quad {\pier {} 
	+ \frac{1}{2} \int_{0}^{t} 
	\left| \frac{\mbox{\boldmath $ f$} (s+h)
	-\mbox{\boldmath $ f$} (s)}{h} 
	\right|_{\mbox{\scriptsize \boldmath $ H$}}^2 ds }
	+c_{\tilde{\delta }} \int_{0}^{t} 
	\left| \frac{\mbox{\boldmath $ v$}_\varepsilon (s+h)
	-\mbox{\boldmath $ v$}_\varepsilon (s)}{h} 
	\right|_{\mbox{\scriptsize \boldmath $ V$}_0^*}^2
	ds.
\end{align*} 
Now, we need that the first two terms in the right hand side remain bounded, uniformly for 
all $h>0$ sufficiently small. Therefore, we go back to \eqref{ape} and integrate it 
from $0$ to $h$, then test by $(\mbox{\boldmath $ v$}_\varepsilon (h)-\mbox{\boldmath $ v$}_0)/h^2$ getting 
\begin{align}
\lefteqn{ 
	\frac{\varepsilon}{2}  \left| 
	\frac{\mbox{\boldmath $ v$}_\varepsilon (h)-\mbox{\boldmath $ v$}_{0}}{h} 
	\right|_{\mbox{\scriptsize \boldmath $ H$}_0}^2 
	+ \frac{1}{2}  \left| 
	\frac{\mbox{\boldmath $ v$}_\varepsilon (h)-\mbox{\boldmath $ v$}_{0}}{h} 
	\right|_{\mbox{\scriptsize \boldmath $ V$}_0^*}^2 
} \nonumber 
	\\
	& \le - \left( {\pier 
	\frac{\mbox{\boldmath $ v$}_\varepsilon (h)-
	\mbox{\boldmath $ v$}_0}{h}	, \frac{1}{h} \int_{0}^{h} 
    \mbox{\boldmath $ P $}\bigl( \partial \varphi \bigl( 
	\mbox{\boldmath $ v$}_\varepsilon (s) 
	\bigr)
	+ \mbox{\boldmath $ \beta $}_\varepsilon 
	\bigl( 
	\mbox{\boldmath $ u$}_\varepsilon (s) 
	\bigr) 
	+ \mbox{\boldmath $ \pi $}
	\bigl( 
	\mbox{\boldmath $ u$}_\varepsilon (s) 
	\bigr) 
	- \mbox{\boldmath $ f$}(s)
	} \bigr) ds
	\right)_{\! \mbox{\scriptsize \boldmath $ H$}_0} \nonumber \\
	& \le 
	\frac{1}{4}  \left| 
	\frac{\mbox{\boldmath $ v$}_\varepsilon (h)-\mbox{\boldmath $ v$}_{0}}{h} 
	\right|_{\mbox{\scriptsize \boldmath $ V$}_0^*}^2 
	+ {\pier \left|
	 \frac{1}{h} \int_{0}^{h} 	 \bigl( \partial \varphi \bigl( 
	\mbox{\boldmath $ v$}_\varepsilon (s) 
	\bigr)
	+ \mbox{\boldmath $ \beta $}_\varepsilon 
	\bigl( 
	\mbox{\boldmath $ u$}_\varepsilon (s) 
	\bigr) 
	+ \mbox{\boldmath $ \pi $}
	\bigl( 
	\mbox{\boldmath $ u$}_\varepsilon (s) 
	\bigr) 
	- \mbox{\boldmath $ f$}(s)
	 \bigr) ds
	\right| _{\mbox{\scriptsize \boldmath $ V$}}^2, } \label{pier11}
\end{align}
{\pier where we have used some properties of the projection operator $\mbox{\boldmath $ P$}$ defined by \eqref{pier10}, in particular  
that $ \vert \mbox{\boldmath $ P$}\mbox{\boldmath $ z $}\vert_{\mbox{\scriptsize \boldmath $ V$}_0} \leq \vert \mbox{\boldmath $ z $}\vert_{\mbox{\scriptsize \boldmath $ V$}}$ for all $ \mbox{\boldmath $ z $}\in \mbox{\boldmath $ V$} $.}
In order that the last quantity {\pier in \eqref{pier11}} be bounded, we need that 
$\mbox{\boldmath $ v$}_0=\mbox{\boldmath $ u$}_0-m_0\mbox{\boldmath $ 1$} \in D(\partial \varphi )$ and especially that 
\begin{gather*}
	-\partial \varphi (\mbox{\boldmath $ v$}_0)
	-\mbox{\boldmath $ \beta $}_\varepsilon (\mbox{\boldmath $ u$}_0)-
	\mbox{\boldmath $ \pi $}(\mbox{\boldmath $ u$}_0)+\mbox{\boldmath $ f$}(0) 
\end{gather*}
remains bounded in {\pier
$\mbox{\boldmath $ V$}$ for all $\varepsilon \in (0,\varepsilon _0]$}
(compare with the assumption (2.40) in the paper \cite{CGS14} {\pier and see the comments just following (2.40)}). 
Therefore, thanks to {\rm (A8)} and {\rm (A9)}{\pier ,} 
the {G}ronwall inequality implies that {\pier the functions}
\begin{gather*}
	t \mapsto \frac{\mbox{\boldmath $ v$}_\varepsilon (t+h)-\mbox{\boldmath $ v$}_\varepsilon (t)}{h}
\end{gather*}
{\pier are} bounded in $L^\infty (0,T-h;\mbox{\boldmath $ V$}_0^*) 
\cap L^2(0,T-h;\mbox{\boldmath $ V$}_0)$, and  
\begin{gather*}
	t \mapsto \varepsilon ^{1/2} \, \frac{\mbox{\boldmath $ v$}_\varepsilon (t+h)-\mbox{\boldmath $ v$}_\varepsilon (t)}{h}
\end{gather*}
{\pier are} bounded in $L^\infty (0,T-h;\mbox{\boldmath $ H$}_0)$ 
uniformly {\pier with respect to} $\varepsilon \in (0,\varepsilon _0]$, 
so that passing to the limit as $h \to 0$ we obtain the bounds
\begin{gather}
	{\pier \varepsilon ^{1/2}|\mbox{\boldmath $ v$}_\varepsilon'|_{L^\infty(0,T;\mbox{\scriptsize \boldmath $ H$}_0)}} + |\mbox{\boldmath $ v$}_\varepsilon '|_{L^\infty (0,T;\mbox{\scriptsize \boldmath $ V$}_0^*)}
	+ 
	|\mbox{\boldmath $ v$}_\varepsilon '|_{L^2(0,T;\mbox{\scriptsize \boldmath $ V$}_0)}
	\le {\pier M_9},
	\label{add1}
\end{gather}
where ${\pier M_9}$ is a positive constant, independent of $\varepsilon \in (0,\varepsilon _0]$. {\pier
From estimate \eqref{add1} it follows that {\pier the estimates for $\{  \beta _\varepsilon (u_\varepsilon ) \}$ and 
$\{ \beta _{\Gamma, \varepsilon} (u_{\Gamma, \varepsilon} ) \}$ in \eqref{pier6}  
can be improved, since the left hand side of \eqref{ape} is now bounded in  
$L^\infty (0,T;\mbox{\boldmath $ V$}_0^*)$ so that
we can test \eqref{ape} by $\mbox{\boldmath $ v$}_\varepsilon $ (bounded in 
$L^\infty (0,T;\mbox{\boldmath $ V$}_0))$ and argue as in \eqref{pier12}, but without a final integration in time. Proceeding in this way, we obtain}  
\begin{gather*}
	\bigl| \beta _\varepsilon (u_\varepsilon ) 
	\bigr|_{L^\infty (0,T;L^1(\Omega ))}
	+ \bigl| \beta _{\Gamma, \varepsilon} (u_{\Gamma, \varepsilon} ) 
	\bigr|_{L^\infty (0,T;L^1(\Gamma))}  \le {\pier M_{10}}.
\end{gather*}
Also, the estimates \eqref{pier13} can be improved to 
\begin{gather*}
	|\omega _\varepsilon |_{L^\infty (0,T)} + 	\bigl | \mbox{\boldmath $ \mu $}_\varepsilon 
	\bigr |_{L^\infty (0,T;\mbox{\scriptsize \boldmath $ V$} )} \le {\pier M_{11}},
\end{gather*}
since $\mbox{\boldmath $ v$}_\varepsilon '$ is bounded in $L^\infty (0,T;\mbox{\boldmath $ V$}_0^*)$ and $\mbox{\boldmath $ f$} \in L^\infty (0,T;\mbox{\boldmath $ H$})$.}
At this point, we can proceed analogously in modifying Lemmas~4.4 and 4.5
{\pier without making the final integrations over $(0,T)$ in the proofs, but 
deducing $L^\infty $ bounds. In particular we find that} 
\begin{gather*}
	\bigl| \beta _\varepsilon (u_\varepsilon) \bigr |_{L^\infty (0,T;H)} 
	\le {\pier M_{12}},
	\\
	|\partial_\nu v_\varepsilon |_{L^\infty (0,T;H_\Gamma )} 
	\le {\pier M_{12}},
	\\
	\bigl |\beta _{\Gamma, \varepsilon} (u_{\Gamma,\varepsilon }) 
	\bigr |_{L^\infty (0,T;H_\Gamma )} \le {\pier M_{12}},
	\\
	|\mbox{\boldmath $ v$}_{\varepsilon} |_{L^\infty (0,T;\mbox{\scriptsize \boldmath $ W$})}
	\le {\pier M_{12}}
\end{gather*}
{\pier in the specified order. Of course, ${\pier M_{10}}$, ${\pier M_{11}}$ and ${\pier M_{12}}$ are positive constants independent of 
$\varepsilon \in (0,\varepsilon _0]$.} 
Now, in the subsequent passage to the limit as $\varepsilon \to 0$, 
for the solution we derive the additional properties that 
\begin{align*}
	&\mbox{\boldmath $ v$}  \in W^{1,\infty }(0,T;\mbox{\boldmath $ V$}^*_0) 
	\cap H^1(0,T;\mbox{\boldmath $ V$}_0) \cap L^\infty (0,T;\mbox{\boldmath $ W$}),
	\\
	&\mbox{\boldmath $ \mu $} \in L^\infty (0,T;\mbox{\boldmath $ V$}),
	\quad \
	{\pier \mbox{\boldmath $ \xi $} \in L^\infty (0,T;\mbox{\boldmath $ H$}).}
\end{align*}
Under these regularities, \eqref{d1} can be rewritten as (cf.~{\pier Remark~2})
\begin{multline}
	\int_{\Omega }^{} {\pier \partial_t v}(t)z dx
	+\int_{\Gamma }^{} {\pier \partial_t v_\Gamma } (t)z_\Gamma  d\Gamma
	+\int_{\Omega }^{} \nabla \mu(t) \cdot \nabla z dx 
	+\int_{\Gamma }^{} \nabla _\Gamma \mu _\Gamma(t) \cdot \nabla _\Gamma z_\Gamma d\Gamma =0
	\\
	\quad \mbox{for all } \mbox{\boldmath $ z$} = (z,z_\Gamma) \in \mbox{\boldmath $ V$}, \ 
	\hbox{for a.a. } t\in (0,T).
	\label{wf1add}
\end{multline} 
Then, taking {\pier a test element $\mbox{\boldmath $ z$} = (z,0) $, with $z \in {\mathcal D}(Q)$, in \eqref{wf1add} and integrating over $(0,T)$, we are led to the equation 
\begin{gather*}
	{\partial_t v} -\Delta \mu = 0 
	\quad \mbox{in } {\mathcal D}'(Q).
\end{gather*}
This implies that $\Delta \mu \in L^2(0,T;H^1(\Omega))$,
due to the regularity of $\partial_t v$.  Moreover, thanks to the fact that
$\partial_t u =\partial_t v $, we also have 
$
	{\pier \partial_t u}-\Delta \mu = 0 $
a.e.\ in $ Q$. 
Furthermore, in view of $\mu _\Gamma \in L^\infty (0,T;V_\Gamma )$, 
we obtain the regularities (see, e.g.,\ \cite[Thm.~3.2, p.~1.79]{BG87} 
and~\cite[Thm.~2.25, p.~1.62]{BG87})
\begin{gather*}
	\mu \in L^2 \bigl( 0,T;H^{3/2}(\Omega ) \bigr),\quad 
	\partial_\nu \mu \in L^2(0,T;H_\Gamma )
\end{gather*}
by the trace theorem. 
The last condition implies that the boundary equation
\begin{gather*}  
	{\pier \partial_t u_\Gamma }+\partial_\nu \mu -\Delta _\Gamma \mu _\Gamma =0
	\quad \mbox{ holds true in } L^2 (0,T;H_\Gamma )  \mbox{ and a.e.\ on }\Sigma,
\end{gather*} 
and provides additional regularity for 
$\mu _\Gamma $:
$
	\mu_\Gamma  \in L^2 \bigl( 0,T;H^2(\Gamma) \bigr),
$
which even yields 
\begin{gather}
	\mu \in L^2 \bigl( 0,T;H^{5/2}(\Omega) \bigr).  \label{pier16}
\end{gather}
Thus, all the regularities for the strong solution specified in Definition~2.2 
are finally obtained, with the further regularity \eqref{pier16}, and equations  and conditions~\eqref{d4}--\eqref{d8} are satisfied.} 
\hfill $\Box$

\section{Appendix}
\setcounter{equation}{0}

We use the same notation as in the previous {\pier sections} for function spaces.

\paragraph{Lemma A.} {\it Let $\Omega \subset \mathbb{R}^d$, $d=2,3$, be bounded, connected and {\pier smooth enough}. 
Then, there exists $c_p>0$ such that} 
\begin{gather*}
	c_p |\mbox{\boldmath $ z$}|_{\mbox{\scriptsize \boldmath $ V$}}^2 
	\le |\mbox{\boldmath $ z$}|_{\mbox{\scriptsize \boldmath $ V$}_0}^2
	\quad \mbox{for all }\mbox{\boldmath $ z$} \in \mbox{\boldmath $ V$}_0.
\end{gather*}

\paragraph{Proof.} 
It {\pier suffices} to show that there exist $C_p>0$ such that 
\begin{gather*}
	|z|_{H}^2 + |z_\Gamma |_{H_\Gamma }^2 \le C_p \Bigl( 
	|\nabla z|_{H^d}^2 + |\nabla_\Gamma  z_\Gamma |_{H_\Gamma ^d}^2 \Bigr) 
	\quad \mbox{for all }\mbox{\boldmath $ z$}:=(z,z_\Gamma ) \in \mbox{\boldmath $ V$}_0.
\end{gather*}
If {\pier this is not the case}, for any $n \in \mathbb{N}$ one can find $\tilde{\mbox{\boldmath $ z$}}_n :=(\tilde{z}_n,\tilde{z}_{\Gamma,n})\in \mbox{\boldmath $ V$}_0$ such that 
\begin{gather*}
	|\tilde{z}_n|_{H}^2 + |\tilde{z}_{\Gamma,n}|_{H_\Gamma }^2 > n \Bigl( 
	|\nabla \tilde{z}_n|_{H^d}^2 + |\nabla_\Gamma  \tilde{z}_{\Gamma,n}|_{H_\Gamma ^d}^2 \Bigr). 
\end{gather*}
Taking $z_n
:=\tilde{z}_n/(|\tilde{z}_n|_{H}^2 
+ |\tilde{z}_{\Gamma,n}|_{H_\Gamma }^2)^{1/2}$, then 
$(z_n)_{|_\Gamma }= \tilde{z}_{\Gamma,n}/(|\tilde{z}_n|_{H}^2 
+ |\tilde{z}_{\Gamma,n}|_{H_\Gamma }^2)^{1/2}$. Therefore,  
{\pier if we set} $\mbox{\boldmath $ z$}_n:=(z_n,z_{\Gamma,n})$ 
with $z_{\Gamma,n}:=(z_n)_{| \Gamma }${\pier ,}
we see that $\mbox{\boldmath $ z$}_n \in \mbox{\boldmath $ V$}_0$ and 
\begin{gather}
	1=|z_n|_{H}^2 + |z_{\Gamma,n}|_{H_\Gamma }^2 > n \Bigl( 
	|\nabla z_n|_{H^d}^2 + |\nabla_\Gamma  z_{\Gamma,n}|_{H_\Gamma ^d}^2 \Bigr)
	\quad \mbox{for all } n \in \mathbb{N}{\pier ;}
	\label{ap0}
\end{gather}
{\pier hence}, we have 
\begin{gather*}
	|z_n|_H^2 \le 1, 
	\quad 
	|z_{\Gamma ,n}|_{H_\Gamma }^2 \le 1, 
	\quad 
	|\nabla z_n|_{H^d}^2 \le \frac{1}{n}, 
	\quad 
	|\nabla _\Gamma z_{\Gamma ,n}|_{H_\Gamma^d}^2 \le \frac{1}{n}
	\quad \mbox{for all } n\in \mathbb{N}.
\end{gather*}
These {\pier inequalities imply} that there exist {\pier a subsequence} of $n$ (not relabeled) and some limit function $\mbox{\boldmath $ z$}:=(z,z_\Gamma ) \in \mbox{\boldmath $ V$}$ 
such that 
\begin{gather}
	z_n \to z \quad \mbox{weakly in } V, \ 
	\mbox{strongly in } H, 
	\label{ap1}\\
	z_{\Gamma, n} \to z_\Gamma  \quad \mbox{weakly in } V_\Gamma, \ 
	\mbox{strongly in } H_\Gamma 
	\label{ap2}
\end{gather}
{\pier as $n \to \infty$,} 	
\begin{gather} 
	\nabla z = 0 \quad \mbox{a.e.\ in } \Omega, \quad 
	\nabla _\Gamma z_\Gamma =0 \quad \mbox{a.e.\ on } \Gamma, 
	\label{ap3}\\
	\int_{\Omega }^{}z dx + \int_{\Gamma }^{}z_\Gamma d\Gamma =0{\pier ;}
	\label{ap4}
\end{gather}
here the {\pier property \eqref{ap4} implies that}
$\mbox{\boldmath $ z$} \in \mbox{\boldmath $ V$}_0$. 
Now, {\pier since $\Omega$ is connected condition~\eqref{ap3} entails} 
that $z$ and $z_\Gamma$ are {\pier the
same constant function}. From \eqref{ap0}--\eqref{ap2} we have 
$|\mbox{\boldmath $ z$}|_{\mbox{\scriptsize \boldmath $ H$}_0}=1$, that is, 
$z \ne 0$ and $z_\Gamma \ne 0$, however in contradiction with \eqref{ap4}. 
This completes the proof of Lemma A. \hfill $\Box$

\paragraph{Lemma B.} {\it {\pier $\mbox{\boldmath $ V$}_0 $ is densely and compactly embedded into $ \mbox{\boldmath $ H$}_0$.}}

\paragraph{Proof.} The strategy of the proof is essentially same as {{\pier in} \cite[Prop.~A1.1]{CF14}. 
For a fixed $\mbox{\boldmath $ z$}=(z,z_\Gamma ) \in \mbox{\boldmath $ H$}_0$ 
and for $n \in \mathbb{N}$, consider the following elliptic problem:
\begin{gather}
	z_n-\frac{1}{n}\Delta z_n =z 
	\quad \mbox{a.e.\ in } \Omega, 
	\label{b1}
\\
	\frac{1}{n} \partial_\nu z_n + (z_n)_{|_\Gamma} =z_\Gamma 
	\quad \mbox{a.e.\ on } \Gamma.
	\label{b2}
\end{gather}
{\pier Is is not difficult to check that (see, e.g.,  \cite[Prop.~A1.1]{CF14}) this sequence 
$\mbox{\boldmath $ z$}_n:=(z_n,z_{\Gamma,n})$, with $z_{\Gamma,n} := (z_n)_{|_\Gamma }$, satisfies}
\begin{gather*}
	\mbox{\boldmath $ z$}_n \to \mbox{\boldmath $ z$} 
	\quad \mbox{strongly in } \mbox{\boldmath $ H$}
	\quad \mbox{as } n \to +\infty. 
\end{gather*}
Moreover, thanks to $\mbox{\boldmath $ z$} \in \mbox{\boldmath $ H$}_0$ we obtain 
\begin{align*}
	\int_{\Omega }^{}z_n dx + \int_{\Gamma }^{} z_{\Gamma,n} d\Gamma 
	& =
	\int_{\Omega }^{} \left( z + \frac{1}{n} \Delta z_n \right) dx 
	+ \int_{\Gamma }^{} \left( z_\Gamma {\pier {} -{}} \frac{1}{n} \partial_\nu z_n \right) d\Gamma 
	\\
	& = \left( \int_{\Omega }^{} z dx + \int_{\Gamma }^{} z_\Gamma d\Gamma \right) 
	+ \frac{1}{n} \left( \int_{\Omega}^{} {\rm div} \nabla z_n dx {\pier {} - {}}\int_{\Gamma }^{} \partial_\nu z_n d\Gamma \right) \  = \ 0.  
\end{align*}
This means that $\{ \mbox{\boldmath $ z$}_n \}_{n \in \mathbb{N}} 
\subset \mbox{\boldmath $ H$}_0$, that is, 
$\mbox{\boldmath $ V$}_0$ is dense in $\mbox{\boldmath $ H$}_0$. 
Next, from the compact {\pier embeddings} $V \subset H$ and $V_\Gamma \subset H$, 
we easily see that for any bounded sequence 
$\{ \mbox{\boldmath $ z$}_n \}_{n \in \mathbb{N}}$ {\pier in} $\mbox{\boldmath $ V$}_0$, 
the sequences of each component,  $\{ z_n \}_{n \in \mathbb{N}}$ 
and $\{ z_{\Gamma,n } \}_{n \in \mathbb{N}}$, have common 
{\pier subsequences that} converge strongly to {\pier some limit functions $z$ and $z_\Gamma $,}  respectively. Therefore, we complete the proof of Lemma B. \hfill $\Box$

\paragraph{Lemma C.} {\it {\pier Let $\varphi : \mbox{\boldmath $ H$}_0 \to [0,+\infty] $
be defined by \eqref{pier7}. Then the} subdifferential 
$\partial \varphi $ on $\mbox{\boldmath $ H$}_0$ is characterized by} 
$$
{\pier  \partial \varphi (\mbox{\boldmath $ z$})
=(-\Delta z,\partial_\nu z-\Delta _\Gamma z_\Gamma ) 
\quad \mbox{\it with} 
\quad  \mbox{\boldmath $ z$}=(z,z_\Gamma ) \in 
D(\partial \varphi )=\mbox{\boldmath $ W$} \cap \mbox{\boldmath $ V$}_0.}
$$  

\paragraph{Proof.} 
{\pier Let} $\mbox{\boldmath $ z$}^* :=(z^*,z_\Gamma ^*) \in \partial \varphi (\mbox{\boldmath $ z$})$ 
in $\mbox{\boldmath $ H$}_0$. 
Then, from the definition of {\pier subdifferential it is straightforward to} obtain 
\begin{gather}
	(\mbox{\boldmath $ z$}^*, \tilde{\mbox{\boldmath $ z$}}
	)_{\mbox{\scriptsize \boldmath $ H$}_0} 
	= (\nabla z,\nabla \tilde{z})_{H^d} 
	+ (\nabla _\Gamma z_\Gamma ,\nabla _\Gamma \tilde{z}_\Gamma )_{H_\Gamma ^d}
	\quad \mbox{for all } \tilde{\mbox{\boldmath $ z$}}:=(\tilde{z},\tilde{z}_\Gamma ) 
	\in \mbox{\boldmath $ V$}_0. \label{c1}
\end{gather}
For each $\tilde{z} \in {\mathcal D}(\Omega )$, 
put $\tilde{\mbox{\boldmath $ z$}}=(\tilde{z},0)$ and 
{\pier take} the test function $\mbox{\boldmath $ P$} \tilde{\mbox{\boldmath $ z$}} 
\in \mbox{\boldmath $ V$}_0$ {\pier above. Using} the property \eqref{proj} we infer that} 
\begin{align*}
	\int_{\Omega}^{} z^* \tilde{z} dx 
	& = 
	(\mbox{\boldmath $ z$}^*, \tilde{\mbox{\boldmath $ z$}}
	)_{\mbox{\scriptsize \boldmath $ H$}} {\pier {}= 
	(\mbox{\boldmath $ z$}^*, \mbox{\boldmath $ P$} \tilde{\mbox{\boldmath $ z$}}
	)_{\mbox{\scriptsize \boldmath $ H$}_0} }\\
	& = \bigl( 
	\nabla z,\nabla \bigl( \tilde{z}-m(\tilde{\mbox{\boldmath $ z$}}) \bigr) \bigr)_{H^d} 
	+ \bigl( \nabla _\Gamma z_\Gamma ,\nabla _\Gamma  \bigl( 0-m(\tilde{\mbox{\boldmath $ z$}}) \bigr) \bigr)_{H_\Gamma ^d} {\pier {}= (\nabla z,\nabla \tilde{z})_{H^d}} {\pier ,}
\end{align*}
namely, $z^* = -\Delta z$ in ${\mathcal D}'(\Omega )$. {\pier Now, as $\mbox{\boldmath $ z$}= (z,z_\Gamma ) \in D(\varphi)=\mbox{\boldmath $ V$}_0$, we have that the trace 
$z_{|_\Gamma} = z_\Gamma  $ is in  $V_ \Gamma$ and 	$-\Delta z = z^* $ lies in $H$. 
Hence, from the theory for elliptic regularity (see, e.g.,\ \cite[Thm.~3.2, p.~1.79]{BG87}), it  follows that $z \in H^{3/2}(\Omega )$. In turn, the trace theory implies that 
$\partial_\nu z \in H_\Gamma$ (cf.~\cite[Thm.~2.25, p.~1.62]{BG87}). Then, recalling~\eqref{c1}, for a general test function  $\tilde{\mbox{\boldmath $ z$}} \in \mbox{\boldmath $ V$}_0$ we have that 
\begin{align*} 
	(\mbox{\boldmath $ z$}^*, 
	\tilde{\mbox{\boldmath $ z$}})_{\mbox{\scriptsize \boldmath $ H$}_0} 
    & =  (z^* , \tilde{z})_H + ( z^*_\Gamma ,\tilde{z}_\Gamma )_{H_\Gamma } \\
	& = (\nabla z,\nabla \tilde{z})_{H^d} 
	+ (\nabla _\Gamma z_\Gamma ,\nabla _\Gamma \tilde{z}_\Gamma )_{H_\Gamma ^d} \\
	& = -(\Delta z, \tilde{z})_H + (\partial_\nu z,\tilde{z}_\Gamma )_{H_\Gamma } 
+ (\nabla _\Gamma z_\Gamma ,\nabla _\Gamma \tilde{z}_\Gamma )_{H_\Gamma ^d},
\end{align*} 
whence, being 	$-\Delta z = z^* $ in $H$, we deduce that
\begin{align*} 
  ( z^*_\Gamma - \partial_\nu z ,\tilde{z}_\Gamma )_{H_\Gamma } = (\nabla _\Gamma z_\Gamma ,\nabla _\Gamma \tilde{z}_\Gamma )_{H_\Gamma ^d}.
\end{align*} 
The last equality implies that $-\Delta _\Gamma z_\Gamma = z_\Gamma ^* -  \partial_\nu z$
is in $H_\Gamma $  and consequently $z_\Gamma \in H^2(\Gamma)$ follows from the boundary 
version of the elliptic regularity theory. Since $z_\Gamma $ is the trace of $z$ on 
the boundary $\Gamma$ and is {\revis sufficiently} smooth (indeed, $z_\Gamma \in H^{3/2}(\Gamma)$), by the elliptic regularity again we conclude that $z \in H^2(\Omega )$, which finally leads to 
$D(\partial \varphi )=\mbox{\boldmath $ W$}\cap \mbox{\boldmath $ V$}_0$.} \hfill $\Box$

\section*{Acknowledgments}

{\pier The} authors wish to express their heartfelt gratitude to 
professors {G}oro {A}kagi and {U}lisse {S}tefanelli, 
who kindly gave them {\pier the opportunity of exchange visits}
supported by the JSPS--CNR bilateral joint research 
\emph{{I}nnovative {V}ariational {M}ethods for {E}volution {E}quations}.
{\pier The present note also benefits from a partial support of the MIUR--PRIN 
Grant 2010A2TFX2 ``Calculus of variations'' and the GNAMPA (Gruppo 
Nazionale per l'Analisi Matematica, la Probabilit\`a e le loro Applicazioni) 
of INdAM (Istituto Nazionale di Alta Matematica) for PC.}


\end{document}